\documentclass[11pt,pra,aps,epsfig,psfig,multicols,showpacs,tightenlines,floatfix]{revtex4}
\usepackage{float}
\usepackage{graphics,bm}
\usepackage{graphicx}
\usepackage{amsmath, amssymb, graphics}
\allowdisplaybreaks

\newcommand{\beq}{\begin{equation}}
\newcommand{\eeq}{\end{equation}}
\newcommand{\bqa}{\begin{eqnarray}}
\newcommand{\eqa}{\end{eqnarray}}
\def\gsim{\mathrel {\vcenter {\baselineskip 0pt \kern 0pt
\hbox{$>$} \kern 0pt \hbox{$\sim$} }}}
\def\lsim{\mathrel {\vcenter {\baselineskip 0pt \kern 0pt
\hbox{$<$} \kern 0pt \hbox{$\sim$} }}}

\protect

\begin{document}

 \title{High accuracy power series method for solving scalar, vector,
 and inhomogeneous nonlinear Schr\"odinger equations}


\author{L. Al Sakkaf  and U. Al Khawaja}

\address{
\it Department of Physics, United Arab Emirates University, P.O.
Box 15551, Al-Ain, United Arab Emirates}


\begin{abstract} \noindent
We develop a high accuracy power series method for solving partial
differential equations with emphasis on the nonlinear Schr\"odinger
equations. The accuracy and computing speed can be systematically and arbitrarily
increased to orders of magnitude larger than those of other methods. 
Machine precision accuracy can be easily reached and sustained for long evolution times within rather short computing time. In-depth analysis and
characterisation for all sources of error are performed by comparing the
numerical solutions with the exact analytical ones. Exact and approximate boundary conditions are considered and shown to minimise errors for solutions with finite background. The method is
extended to cases with external potentials and coupled nonlinear
Schr\"odinger equations.
 \end{abstract}

\maketitle

\section{Introduction} \label{intsec}
The nonlinear Schr\"odinger equation (NLSE) is truly a universal
equation as it describes major fields such as Bose-Eienstein condensation \cite{hcbook},
nonlinear optics \cite{opticsbooks}, ocean waves \cite{khareef}, and
many others \cite{others1,others2,others3}. This has stimulated
extensive interest in its analytical \cite{ourbook} and numerical
solutions
\cite{meth2,meth3,meth12,meth4,meth6,meth8,meth9,meth11,meth7,meth10,meth5}.
Over decades, knowledge about its analytical solutions has
accumulated such that it is now rare to find a new solution
\cite{ourbook}. A plethora of numerical methods have also been developed to
solve its nonlinear initial value problem, such as the incoherent
scattering of solitons with each other or the scattering of solitons by external
potentials.  Some solutions demand higher accuracy than others such
as dark solitons or vortex excitations since they have nonzero
background or extend over the whole system, which warrants accurate
account of the boundary conditions. Some other solutions have a fast
time evolution or high curvature such as the coalescing soliton
molecule or Peregrine soliton. This kind of solutions requires
accurate integrator of the time derivative. Many of the numerical
methods developed already solve these problems, but there is always
a demand on increasing accuracy and decreasing computing time and
memory cost, particularly for large system sizes or long evolution
times. Explicit methods solving the NLSE can be categorised into two
major classes, namely spectral methods and finite difference methods
\cite{meth12}, in addition to other methods using, for instance,
quadrature discretisation \cite{32,33,34,35} or wavelet expansion
\cite{meth7}.

Here, we present a method that can systematically increase the
accuracy in both the spacial and temporal axes. For the temporal
evolution, we use an iterative power series method that we have
developed previously for ordinary differential equations
\cite{pspaper} and applied later to fluid flow \cite{fluid}. The
accuracy in the time evolution increases with the maximum power in the time power series, $s$. For
the spacial part, we use a $p$-point stencil to discretise the
second derivative, where $p\ge3$ is a positive odd  integer. The accuracy
can be systematically increased by increasing $s$ and $p$. The
method requires the knowledge of the initial profile and its
boundary conditions. While an arbitrary initial profile can be used,
using an exact solution as the initial profile, makes it possible to
calculate the evolution of error and compare it with other numerical
methods. The exact solutions we consider here include: moving bright
soliton, moving dark soliton, Peregrine soliton, and soliton
molecule. For all of these solutions, our method shows a
remarkable performance with accuracy that can reach the machine
precision for a long evolution time in a rather short computing
time.

We perform an extensive analysis of the different sources of error
originating from spacial discretisation, temporal discretisation,
and boundary conditions. We then compare our method with two methods
representing the finite difference and spectral methods. Among the
many finite difference methods, the so-called  generalized finite-difference time-domain (G-FDTD) method
\cite{phd,meth2,meth3,meth12} is superior in
its high accuracy. Based on our analysis and understanding of the interplay
between the different sources of error, the present method makes
significant enhancements on the G-FDTD method in terms of accuracy
or computing time. Specifically, we enhance on the time stepping
method, make no approximations in the derivation of the recursion
relations of the time power series, and account for the evolution of
boundary points exactly, even when the initial profile is not an
exact solution. Furthermore, we show that using the so-called constant wave (CW) exact solution to compute approximate boundary
conditions is in most cases as good as using the exact ones for localised solutions over a uniform background, as long
as the localisation does not come close to the boundaries
within the considered time domain. This is shown with a detailed
comparison performed for the bright and dark soliton solutions. For
the spectral methods, we compare with the Fourier split-step (SS) method \cite{ss},
where we show that the present method is significantly more accurate
and faster.

The method is extended to inhomogeneous NLSE and applied to the 
nonintegrable case of a bright soliton
scattered by a reflectionless potential well. An accurate accounte
to the quantum reflection effect and to the value of the
critical speed \cite{coodman,brand,brand2} are obtained, where extremely high accuracy
is required when the soliton speed is close enough to the critical
speed. We show that other numerical methods lead to the wrong
outcome (reflection instead of transmission), while the present
method captures the correct behaviour right at its lowest level of
accuracy ($p=3$).

Finally, the method is generalised to the case  of two coupled
NLSEs, known as the Manakov system, from which the  evolution of the
dark-bright soliton is calculated accurately.

The rest of the paper is organised as follows. In Section
\ref{numsec}, we present the proposed theoretical framework and
algorithm of the method. In Section \ref{errsec}, we perform a
detailed analysis of the different sources of error and characterize
them in terms of $s$ and $p$. In Section \ref{compsec}, we compare the
accuracy and CPU run time with other methods. In Section
\ref{othersec}, we consider the Peregrine soliton and soliton molecule as
initial profiles. In Section \ref{vsec}, we extend the method to
NLSE with an external potential. In Section \ref{twosec}, we
generalize to the two-coupled NLSE or Manakov system. We end in
Section \ref{concsec} with a summary and outlook for future work.

\section{Numerical method}
\label{numsec}
While emphasis will be on the NLSE, the method we describe below can be modified to solve other evolution equations. The fundamental NLSE can be written in dimensionless form as
\begin{equation}\label{nlse}
i\, \frac{\partial}{\partial
t}\psi(x,t)+g_1\,\frac{\partial^2}{\partial
x^2}\psi(x,t)+g_2\,|\psi(x,t)|^2\,\psi(x,t)=0,
\end{equation}
where $\psi(x,t)$ is a complex function, $g_1$ and $g_2$ are
arbitrary real constants representing the strength of dispersion and
nonlinear terms, respectively. In nonlinear optics, the NLSE
describes the propagation of pulses in nonlinear media. In such a
context,  the dispersion term corresponds to the  group velocity
dispersion (GVD), which, depending on the sign of $g_1$, compresses
or spreads out the pulse, while the nonlinear term corresponds to
what is known as  the Kerr effect, which describes the modulation of
the refractive index of the medium as a response to the propagating
light pulse.

\begin{figure}[h]
    \centering
    \includegraphics[scale=0.25]{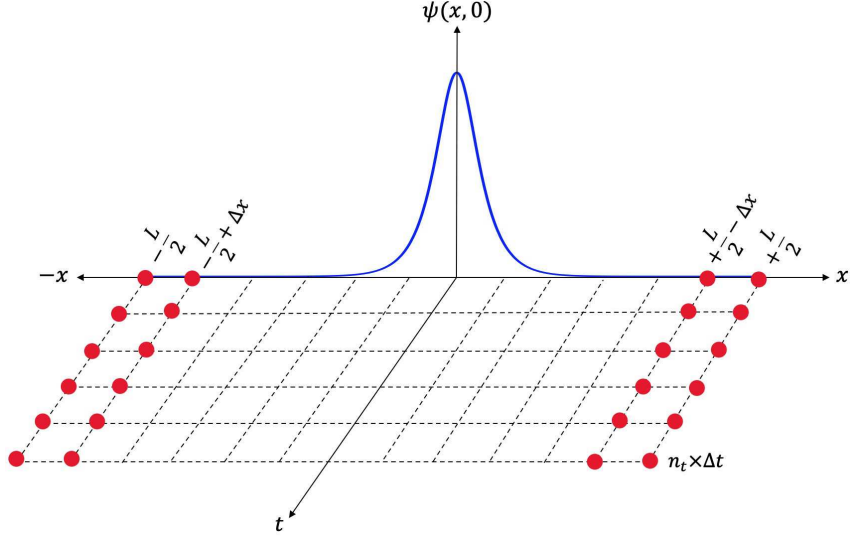}
    \caption{Schematic diagram showing the initial profile (blue curve) and boundary points (red dots) for the case with $p=5$. Number of the boundary points in each side of the grid equals $(p-1)/2$. Each point is a part of a sequence of $n_t$ points along the $t$-axis.}
    \label{fig1}
\end{figure}

The statement of the  problem is defined as follows: Given an
arbitrary initial profile $\psi_0(x)=u_0(x)+i\,v_0(x)$ and boundary
conditions on $\psi(\pm L/2,t)$ at the edges of the spacial domain,
$x=\pm L/2$, what is the time evolution of $\psi_0(x)$ governed by
the NLSE, Eq. (\ref{nlse})? The method presented here solves this
nonlinear initial value problem, which is described briefly as
follows. The solution is expanded in a power series in time as
$\psi(x,t)=c_0(x)+c_1(x)t+c_2(x)t^2+\dots+c_s(x)t^s$, where $s$ is a
positive integer. Recursion relations for the coefficients
$c_l(x),\,\,l>0$, will be given in terms of the initial profile
$c_0(x)$ upon substituting in (\ref{nlse}). The spacial domain is
discretised using a $p$-point stencil to replace the second
derivative, where $p$ is an odd integer $\ge3$. As a result, the
first and last $(p-1)/2$ points of the spacial grid can not be
determined by the recursion relations and need to be determined from
boundary conditions. For the class of solutions which we consider
here, namely a localised profile over a uniform background, the
CW exact solution may be used to
accurately calculate these boundary conditions.  The schematic
figure, Fig. \ref{fig1}, depicts the picture just described.
Accuracy in the method is thus determined by three factors:  i) the
order of the time power series $s$, with error $\propto \Delta
t^{s+1}$, where $\Delta t$ is the small discretisation in the time
domain,  ii) the number of points in the $p$-point formula
approximating the second spacial derivative with error $\propto
\Delta x^{p-1}$, where $\Delta x$ is the small discretisation in the spatial
domain, iii) the accuracy in the boundary conditions, namely
how accurately does the CW solution represent the evolution of the
first and last $(p-1)/2$ boundary points.

We use in this work values of $s\le4$ and mostly $p\le23$. It turns
out that very high accuracy which can easily reach machine precision
is accessible for a long evolution time but still with a short run
time. In the following, the
method is described in detail.

\subsection{Time evolution and recursion relations}
Without loss of generality, we write the general solution in the cartesian complex form of
\begin{eqnarray}
\psi(x,t)&=&u(x,t)+i\,v(x,t),
\end{eqnarray}
 where $u(x,t)$ and $v(x,t)$ being real functions. Inserting in (\ref{nlse}),
 generates the following two equations from the real and imaginary parts
\begin{eqnarray}
g_1\,\frac{\partial^2}{\partial x^2}u(x,t)+g_2\,
\left[u^2(x,t)+v^2(x,t)\right]u(x,t)-\frac{\partial}{\partial t}v(x,t)&=&0,\label{u}\\
g_1\,\frac{\partial^2}{\partial
x^2}v(x,t)+g_2\,\left[u^2(x,t)+v^2(x,t)\right]v(x,t)+\frac{\partial}{\partial
t}u(x,t)&=&0.\label{v}
\end{eqnarray}
To integrate the time operator, we expand the solutions $u(x,t)$ and
$v(x,t)$ in powers of the small discretization of the time domain,
$\Delta t$, as follows
\begin{eqnarray}
u(x,t)&=&\sum_{l=0}^{s}a_l(x)\,{\Delta t}^{\,l},\label{seru}\\
v(x,t)&=&\sum_{l=0}^{s}b_l(x)\,{\Delta t}^{\,l}\label{serv},
\end{eqnarray}
where $s$ is a positive nonzero integer and the function
coefficients are defined by
\begin{eqnarray}
a_l(x)&=&\frac{1}{l!}\,\frac{\partial^l}{\partial t^l}u(x,t)|_{t=0}\label{uvexp1},\\
b_l(x)&=&\frac{1}{l!}\,\frac{\partial^l}{\partial t^l}v(x,t)|_{t=0}\label{uvexp2}.
\end{eqnarray}
The initial profile is given by $\psi_0(x)=a_0(x)+i\,v_0(x)$. We
refer to the maximum power of this expansion, $s$, as the ``order".
The order is one of two factors that can be used to increase the accuracy of the
method. By terminating the power series at $s$, an error
\begin{eqnarray}
error_{s}&=&\mathcal{O}\left(\frac{\Delta
t^{\,s+1}}{(s+1)!}\,\frac{\partial^{s+1}}{\partial
t^{s+1}}u(x,t)|_{t=0}\right)\nonumber\\&=&\mathcal{O}\left(\frac{\Delta
t^{\,s+1}}{(s+1)!}\,\frac{\partial^{s+1}}{\partial
t^{s+1}}v(x,t)|_{t=0}\right)
\label{errors}
\end{eqnarray}
is introduced. To proceed, we need to specify the value of $s$, which
we choose as $s=4$. The expansions (\ref{seru}) and (\ref{serv})
then read
\begin{eqnarray}
u(x,t)&=&a_0(x)+a_1(x)\,t+a_2(x)\,t^2+a_3(x)\,t^3+a_4(x)\,t^4\label{uvexp1},\\
v(x,t)&=&b_0(x)+b_1(x)\,t+b_2(x)\,t^2+b_3(x)\,t^3+b_4(x)\,t^4\label{uvexp5}.
\end{eqnarray}
Substituting back into (\ref{u}) and (\ref{v}), recursion relations
are obtained by equating the coefficients of $\Delta t^{\,i}$,
$i=0,\,1,\,2,\,3$ to zero, as follows

\begin{eqnarray}
\label{a1eq}
a_1   &=&-g_2\,a_0^2   \,b_0   -g_2\,b_0^3   -g_1\,b_0^{\prime\prime},\\
a_2   &=&\frac{1}{2}\left[-2\,g_2\,a_0   \,a_1   \,b_0   -g_2\,a_0^2   \,b_1   -3\,g_2\,b_0^2   \,b_1   -g_1\,b_1^{\prime\prime}   \right]\label{a2eq},\\
a_3   &=&\frac{1}{3}\Big[-g_2\,a_1^2   \,b_0   -2\,g_2\,a_0   \,a_2   \,b_0   -2\,g_2\,a_0   \,a_1   \,b_1   -3\,g_2\,b_0   \,b_1^2   \nonumber\\&&-g_2\,a_0^2   \,b_2   -3\,g_2\,b_0^2   \,b_2   -g_1\,b_2^{\prime\prime}   \Big]\label{a3eq},\\
a_4   &=&\frac{1}{4}\Big[-g_2\,a_1^2   \,b_1   -g_2\,b_1^3   -2\,g_2\,a_0   \,a_3   \,b_0   -2\,g_2\,a_0   \,a_2   \,b_1   -6\,g_2\,b_0   \,b_1   \,b_2   \nonumber\\&&-2\,g_2\,a_1   \,a_2   \,b_0   -2\,g_1\,a_1   \,a_0   \,b_2   -g_2\,a_0^2   \,b_3   -3\,g_2\,b_0^2   \,b_3   -g_1b_3^{\prime\prime}   \Big]\label{a4eq},\\
b_1   &=&g_2\,a_0   \,b_0^2   +g_2\,a_0^3   +g_1\,a_0^{\prime\prime},
\\
b_2   &=&\frac{1}{2}\left[2\,g_2\,a_0   \,b_0   \,b_1   +g_2\,b_0^2   \,a_1   +3\,g_2\,a_0^2   \,a_1   +g_1\,a_1^{\prime\prime}   \right]\label{b2eq},
\\
b_3   &=&\frac{1}{3}\Big[g_2\,a_0   \,b_1^2   +2\,g_2\,a_0   \,b_2   \,b_0   +2\,g_2\,a_1   \,b_0   \,b_1   +3\,g_2\,a_0   \,a_1^2   \nonumber\\&&+g_2\,a_2   \,b_0^2   +3\,g_2\,a_0^2   \,a_2   +g_1\,a_2^{\prime\prime}   \Big],\\
b_4   &=&\frac{1}{4}\Big[g_2\,a_1   \,b_1^2   +g_2\,a_1^3   +2\,g_2\,a_0   \,b_3   \,b_0   +2\,g_2\,a_0   \,b_1   \,b_2   +6\,g_2\,a_0   \,a_1   \,a_2   \nonumber\\&&+2\,g_2\,a_1   \,b_2   \,b_0   +2\,g_2\,a_2   \,b_0   \,b_1   +g_2\,b_0^2   \,a_3   +3\,g_2\,a_0^2   \,a_3   +g_1a_3^{\prime\prime}   \Big]\label{b4eq},
\end{eqnarray}
where we hid the $x$-dependence for convenience and
$(\cdot)^{\prime\prime}$ denotes a second derivative with respect to
$x$. The structure of the recursion relations is such that the
$l^{\rm th}$  coefficient is determined by the $(l-1)^{\rm th}$
coefficients and the second derivative of one of the $(l-1)^{\rm
th}$ coefficients. For instance, $a_1$ is given in terms of $a_0$,
$b_0$, and $b_0^{\prime\prime}$.

\subsection{Spacial discretisation and boundary conditions}
The spacial domain of size $[-L/2,L/2]$ is divided into a number
$n_x$ of spacial discretizations, $\Delta x=L/(n_x-1)$ such that
$x_i=-L/2+(i-1)\times\Delta x,\,\,i=1,2,\dots,n_x$. The coefficients
$a_{1-4}(x)$ and $b_{1-4}(x)$ are then discretized and denoted as
$a_{1-4}^i$ and $b_{1-4}^i$, respectively. The crucial point is how
to discretize the second derivative of the coefficients. The lowest
order discretization is given by the three-central-point (Euler)
formula $a_{1-3}^{\prime\prime}(x)=\left[a_{1-3}(x+\Delta
x)+a_{1-3}(x-\Delta x)-2a_{1-3}(x)\right]/\Delta x^{\,2}$, which in discretized form reads
$a_{1-3}^{\prime\prime}=(a_{1-3}^{i+1}+a_{1-3}^{i-1}-2a_{1-3}^i)/\Delta
x^{2}$, and similarly for $b_{1-3}^{\prime\prime}$. Using the forward or backward formula would be as accurate
as the central point formula, but the later is more convenient for
clearly exhibiting the symmetry in the boundary conditions, as will
be detailed below.

A major source of error is introduced by the discretization of the
second derivative. For instance, in deriving the three-point
formula, an error of order $\Delta x^2$ is introduced. For higher
accuracy, we use the following formula for a general central point
second derivative, denoted here and throughout as
$f_p^{\prime\prime}$, of any odd number of points, $p$,
\begin{eqnarray}
f_p^{\prime\prime}\equiv\frac{d^2f(x)}{dx^2}&=&\frac{1}{\Delta
x^2\sum_{j=1}^{(p-1)/2}C_j{j^2}}\nonumber\\&&\times\sum_{j=1}^{(p-1)/2}C_j\left[f(x+j\Delta
x)+f(x-j\Delta
x)-2f(x)\right]\nonumber\\&&+\mathcal{O}\left(\frac{(\Delta x)^{p-1}}{(p+1)!}\frac{df^{p+1}(x)}{dx^{p+1}}\right)\label{eqapp3p},
\end{eqnarray} where the $p$-dependent
coefficients, $C_j$, are determined by an appropriate linear
combination of the Taylor expansions of $f(x+j\,\Delta x),\,
j=\pm1,\,\pm2,\,\dots,\,\pm(p-1)/2 $, which results in the following system
\begin{equation}
\sum_{j=1}^{(p-1)/2}C_j{j^{2i}}=0\label{eqapp2}
\end{equation}
for $i=2,3,\dots,(p-1)/2$, to be solved in terms of $C_j$. The last
term in (\ref{eqapp3p}) gives the order of error introduced by this
approximation to the second derivative.
This error appears in the coefficients $a_{1-4}$ and $b_{1-4}$ of
Eqs. (\ref{seru}) and (\ref{serv}). However,  the dominant
contribution will be from $a_1$ and $b_1$ since they are multiplied
by $\Delta t$, while the rest of coefficients are multiplied by
higher powers of $\Delta t$. Therefore,  the second source of error
in our method takes the form
\begin{eqnarray}
error_p&=&\max\left[\Delta t\,\frac{\partial}{\partial t}\frac{(\Delta x)^{p-1}}
{(p+1)!}\frac{\partial |\psi^{p+1}(x,t)|}{\partial x^{p+1}}\right].\label{errorp}
\end{eqnarray}
The detailed derivation of the $p$-point formulae is relegated to Appendix \ref{appa} together with explicit formulae for the cases $p=5,\dots,23$.

Using the $p$-point approximation to the second derivative, the
coefficients $a_{1-4}^i$ and $b_{1-4}^i$ will be given in terms of a
$p$-point stencil composed of $(p-1)/2$ points to the left and
$(p-1)/2$ points to the right of the central point. Consequently,
the recursion relations can only be used for points $x_{(p-1)/2}<
x_i< x_{n_x-(p-1)/2}$, which excludes the first and last $(p-1)/2$
points from the grid. The evolution of these two sets of points
needs to be determined by the boundary conditions. In practise,
there are different possibilities to consider. In many situations,
an exact solution of the fundamental NLSE, (\ref{nlse}), is used as
an initial profile such as the bright or dark soliton to be
scattered by other solitons or by a potential. In this case, the
initial profile is not an exact solution anymore but represents an
accurate approximation near the edges of the spacial domain.  For
these situations, the boundary conditions can be calculated from
(\ref{uvexp1}) and (\ref{uvexp5}) using $u(x,t)$ and $v(x,t)$ as the
exact solution considered for the initial profile. In other
situations, an arbitrary initial localised profile that is not an
exact solution to (\ref{nlse}) is used. Near the edges of the
spacial domain, which are assumed to be sufficiently far from the
localisation for all times, the solution decays to a uniform
background.  The uniform background, namely the CW solution, is an
exact solution to the fundamental NLSE, (\ref{nlse}). Therefore, the
evolution of the boundary points may be calculated from
(\ref{uvexp1}) and (\ref{uvexp5}) using $u(x,t)$ and $v(x,t)$ as the
CW solution. The advantage of using the CW solution over using exact
solutions to calculate the evolution of the boundary points is that
it applies to all localised initial profiles as long as the
localisation is away from the edges. This is very practical since it
means that we can evolve any initial profile without a priori
knowledge of the full time evolution of the boundary points. Using
an exact solution, on the other hand, has the advantage of analytic
continuation allowing the localisation to cross the boundaries of
the spacial domain. We have verified that for cases with evolution
times long enough to perform realistic numerical experiments, using
the CW leads to almost the same accuracy in the evolved profile as
when the exact solution is used for the initial profile.

Based on the above, the evolution of the coefficients $a_l(x)$ and $b_l(x)$ for for the bulk of the grid, namely $1+(p-1)/2\le i\le n_x-(p-1)/2$, read in the discretised form
\begin{eqnarray}
a_1^i   &=&-g_2\left((a_0^i)^2+(b_0^i)^2\right) \,b_0^i   -\frac{g_1}{\Delta x^2}\sum_{j=1}^{(p-1)/2}{\bar C}_j\,\left[b_0^{i+j}+b_0^{i-j}-2b_0^{i}\right]\label{req1},\\
b_1^i   &=&g_2\left((a_0^i)^2 +(b_0^i)^2\right) \,a_0^i   +\frac{g_1}{\Delta x^2}\sum_{j=1}^{(p-1)/2}{\bar C}_j\,\left[a_0^{i+j}+a_0^{i-j}-2a_0^{i}\right],
\label{req2}
\end{eqnarray}
where ${\bar C}_j=C_j/\sum_{j=1}^{(p-1)/2}C_jj^2$. Equations for the
rest of coefficients, $a_{2-4}^i$ and $b_{2-4}^i$, can similarly be
obtained by discretising (\ref{a2eq}-\ref{a4eq}) and
(\ref{b2eq}-\ref{b4eq}). The boundary points, $1\le i\le (p-1)/2$
and $n_x-(p-1)/2< i\le n_x$, are calculated from the boundary
condition as
\begin{eqnarray}
a_l^i&=&\frac{1}{l!}\frac{\partial^l}{\partial t^l}u(x_i,t)|_{t=0}\label{req1b},\\
b_l^i&=&\frac{1}{l!}\frac{\partial^l}{\partial t^l}v(x_i,t)|_{t=0}\label{req2b},
\end{eqnarray}
where, $u(x_i,t)$ and $v(x_i,t)$ correspond to an exact solution of
(\ref{nlse}). They may correspond to the time-dependent localized
solution from which the initial profile is used. Alternatively, they
may be the CW solution that the initial profile approaches at the
boundaries. For the bright soliton, the CW solution that describes
the asymptotes at the boundary is zero. Therefore, for this special
case, all coefficients $a_{1-4}^i$ and $b_{1-4}^i$ can be set to
zero at the boundary points defined above. This saves considerably
on memory and CPU time. However, throughout this paper, we did not
use such trivial boundary conditions. We restricted the boundary
conditions either to the exact localized solution or the CW. A
number of examples on localised solutions are considered in the
following section, while the uniform, CW, solution of the NLSE,
(\ref{nlse}) considered here, is given by
\begin{eqnarray}
\psi(x,t)&=&A_0\,e^{i\,\left[\left(g_2\,A_0^2-\frac{k^2}{4\,g_1}\right)t+\frac{k}{2\,g_1}(x-x_0)\right]}\label{bseq},
\end{eqnarray}
from which we define
\begin{eqnarray}
u(x,t)&=&A_0\,\text{cos}\left[\left(g_2\,A_0^2-\frac{k^2}{4\,g_1}\right)t+\frac{k}{2\,g_1}(x-x_0)\right],
\end{eqnarray}
\begin{eqnarray}
v(x,t)&=&A_0\,\text{sin}\left[\left(g_2\,A_0^2-\frac{k^2}{4\,g_1}\right)t+\frac{k}{2\,g_1}(x-x_0)\right],
\end{eqnarray}
where $A_0$, $k$, and $x_0$ being arbitrary real constants.

Finally, the method can be summarized with the following
algorithm:\\\\

\begin{center}{\bf Algorithm:}\end{center}
{\it {\bf1)} Initial profile, for $1\le i\le n_x$: $a_0^i=u_0(x_i)$
and
$b_0^i=v_0(x_i)$.\\\\
{\bf2)} Boundary conditions, for $l>0$, $1\le i\le (p-1)/2$ and $n_x-(p-1)/2< i\le n_x$:\\
\begin{eqnarray}
\nonumber
a_l^i&=&\frac{1}{l!}\frac{\partial^l}{\partial t^l}u(x_i,t)|_{t=0},\\
b_l^i&=&\frac{1}{l!}\frac{\partial^l}{\partial
t^l}v(x_i,t)|_{t=0}\nonumber,
\end{eqnarray}
[Eqs. (\ref{req1b}) and (\ref{req2b})]\\\\
{\bf3)} Recursion relations, for  $l>0$, $(p-1)/2< i\le n_x-(p-1)/2$:\\
\begin{eqnarray}\nonumber a_l^i&=&a_l^i\left(a_0^i,b_0^i,\sum_{j=1}^{(p-1)/2}{\bar
C}_j\,\left[b_{l-1}^{i+j}+b_{l-1}^{i-j}-2b_{l-1}^{i}\right]\right),\\
b_l^i&=&b_l^i\left(a_0^i,b_0^i,\sum_{j=1}^{(p-1)/2}{\bar
C}_j\,\left[a_{l-1}^{i+j}+a_{l-1}^{i-j}-2a_{l-1}^{i}\right]\right)\nonumber.\end{eqnarray}
[Eqs. (\ref{req1}), (\ref{req2}), and similar equations for the rest of coefficients.]\\\\
{\bf4)} Time evolution and update:
\begin{eqnarray}
a_0^i&\leftarrow&\sum_{l=0}^{s}a_l^i\,{\Delta t}^{\,l},\label{seru2}\\
b_0^i&\leftarrow&\sum_{l=0}^{s}b_l^i\,{\Delta t}^{\,l}\label{serv2}.
\end{eqnarray}
{\bf5)} Return to step {\bf3} with the updated values of $a_0^i$ and
$b_0^i$. }\\

\begin{figure}[h]
    \centering
    \includegraphics[scale=0.23]{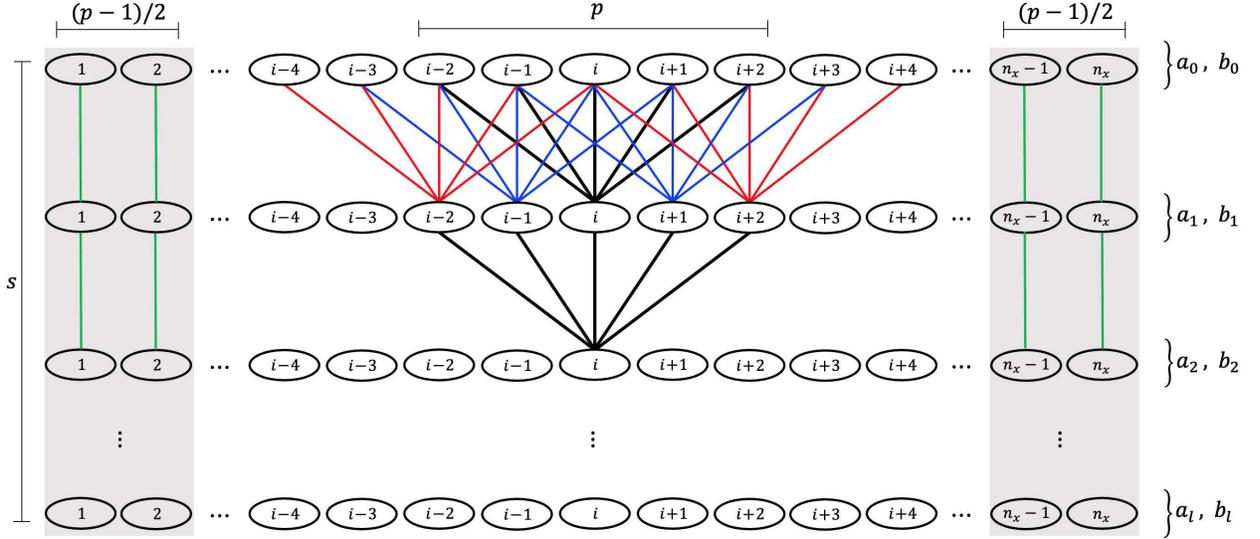}
    \caption{Schematic figure depicting the algorithm
    of the present method. The figure is based on the special case of
    $p=5$. Each coefficient, $a_l$ or $b_l$, is calculated from 5 points of
    the previous order, $a_{l-1}$ and $b_{l-1}$, according to Eqs. (\ref{req1}) and (\ref{req2}),
    and similar equations for the rest of coefficients.
    The first and last two points are calculated by the boundary
    conditions according to Eqs. (\ref{req1b}) and (\ref{req2b}).}
    \label{fig2}
\end{figure}

The algorithm is also depicted schematically in Fig.
\ref{fig2}.

\section{Error analysis}
\label{errsec} The main aim here is to calculate and characterize
the error of our method. As pointed out in the previous section,
there are two main sources of error. The first source of error
arises from the termination of the time power series at order $s$,
namely $error_s$ given by Eq. (\ref{errors}). The second source of
error is due to approximating the second derivative by the $p$-point
formula, $error_p$ given by Eq. (\ref{errorp}). We verify this
understanding through the numerical solution of bright and dark
solitons.

\subsection{Bright soliton}
The exact movable bright soliton solution of  Eq. (\ref{nlse}) can
be expressed as
\begin{eqnarray}\label{bright}
\psi(x,t)&=&A_0\,\sqrt{\frac{2\,g_1}{g_2}}\,\text{sech}\left\{A_0\left[x-(x_0+k\,t)\right]\right\}\,e^{i\left[\frac{k}{2\,g_1}(x-x_0)+\left(\frac{4\,A_0^2\,g_1^2-k^2}{4\,g_1}\right)(t-t_0)+\phi_0\right]}\label{bs},
\end{eqnarray}
where  $g_1 g_2>0$, and the arbitrary real constants  $A_0$, $x_0$, $t_0$, $k$, and $\phi_0$ physically define the height of the wave,  spacial shift, temporal shift, soliton speed, and  global phase, respectively.

The error is defined as
\begin{equation}
error=\max\left[||\psi_{\rm numerical}(x_i,t_f)|-|\psi_{\rm
exact}(x_i,t_f)||\,\right],\hspace{1cm}1\le i\le n_x,
\end{equation}
where $t_f$ is the final time of evolution. In semi-log plots versus
$p$ for the 4 values of $s$, Fig. \ref{fig3} shows the general
behavior of a decreasing error that saturates at a certain value.
Our analysis shows that the decreasing part corresponds to $error_p$
and the saturating part corresponds to $error_s$. For low values of
$p$, the error in the $p$-point formula, $error_p$, is larger than
the error in the order, $error_s$. With larger values of $p$,
accuracy enhances such that $error_p$ becomes smaller than
$error_s$, and thus the total error is dominated by $error_s$ which
is independent of $p$. To verify this understanding, we calculate
$error_p$, as given by Eq. (\ref{errorp}) and plot it with the black
filled circles where it is clear that the theoretical prediction of
this part of error follows the numerical one. The order error,
$error_s$, is calculated from Eq. (\ref{errors}), and is plotted
with the dashed horizontal lines. Here again, the theoretical
prediction for $error_s$ agrees very well with the numerical values.
The figure shows clearly the interplay between the effects of $s$
and $p$ on the accuracy: Decreasing the error with $p$ is limited by
a saturation minimum set by $s$.  We have repeated this calculation
for decreasing time discretization but with keeping the final time
the same. The purpose of this is to verify that the saturation
values do indeed decrease according to  Eq. (\ref{errors}), which is
clearly the case as can be seen in all cases considered. It shows
also, as expected, that decreasing $\Delta t$ has the same effect as
increasing the order. In Fig. \ref{fig4}, we show the effect of
increasing the order on accuracy. The saturation levels correspond
to $error_p$ where error is no longer depending on the order. For
larger $p$, higher order is needed to reach saturation.

For realistic applications, it is important to keep the error small during long times of evolution. We show in Fig. \ref{fig5} the time evolution of error up to $t=40$ with different values of $p$. For $p=3,5,7$, the error grows linearly with time for most of the time interval. For larger values of $p$, the error starts to saturate at a value that decreases with increasing $p$. For $p=23$, the error saturates at the machine precision. Therefore, the numerical solution can be considered as exact up to the machine precision within the time interval considered. On a semi-logarithmic scale, we plot in Fig. \ref{fig6} the error at the end of time evolution, $error(t=t_f)$, versus $p$, which shows how rapidly the error drops to the machine precision with increasing $p$. For longer time evolution, Fig. \ref{fig7}  shows that even with $p=23$, the error starts to grow with time. Larger value of $p$ is needed to get the error back to saturation.

Since using an exact analytical solution to calculate the boundary conditions is not the most general case, we investigate the effect of replacing the exact boundary conditions by approximate ones. For the bright soliton,  all coefficients $a_l^i$ and $b_l^i$ at the boundaries, namely with $1\le i< (p-1)/2$ and $n_x-(p-1)/2< i\le n_x$, are set to zero. In Fig. \ref{fig8}, we show that using approximate boundary conditions leads to an error that is almost identical to that when exact boundary conditions are used.

\subsection{Dark soliton}
Solutions with finite background are typically more demanding computationally due to errors from the edges of the spacial domain. This introduces another source of error. In the present method, boundary points are fixed by boundary conditions through an exact or approximate analytical solutions while the bulk of the spacial grid is evolved according to the numerical method. The difference in evolution procedure  generates high frequency oscillations stemming from the boundary between the points evolved with the numerical method and the points evolved with the boundary conditions.

The dark soliton we consider for comparison is given by
 \begin{equation}\label{dark}
\psi(x,t)=A_0\,\sqrt{\frac{-2\,g_1}{g_2}}\,\text{tanh}\{A_0\,[x-(x_0+k\,t)]\}\,e^{-i\,[-\frac{k}{2\,g_1}\,(x-x_0)+\frac{8\,g_1^2\,A_0^2+k^2}{4\,g_1}\,(t-t_0)+\phi_0]},
\end{equation}
where $g_1\,g_2<0$.
Time evolution of error is shown in Fig. \ref{fig9}. Similar to the bright soliton case, the linear dependence of error on evolution time disappears with increasing $p$. However, it is noticed here that the saturation value for $p=23$ is around $1.5\times 10^{-12}$ which is not at machine precision, as the case was with bright soliton. In Fig. \ref{fig10}, this can also be seen with error at the final evolution time plotted versus $p$. Investigating this behaviour further showed that it is due to the errors at the boundaries. Here, the background is finite and boundary errors appear more prominently unlike the case of zero background for bright soliton. A snapshot of the error is shown in Fig. \ref{fig11} where it is clear that the error is significant only at the boundaries and the centre. The central error is associated with the structure of the dark soliton. The nature of the boundary error is different; it is caused, as mentioned above, by fixing the boundary points to fixed values and evolving the other points using the numerical scheme.  For small $p$, the central error is dominant. Increasing $p$ reduces both the central and boundary errors, but the boundary error saturates after a certain value of $p$, while the central error keeps decreasing.  The boundary error becomes dominant for larger $p$. This is verified in Fig. \ref{fig12} where we plot the central and boundary errors separately. While the central error is responsive to increasing $p$, where it ultimately decreases down to machine precision, the boundary error saturates at a larger value. Nonetheless, for the given parameters which are realistic, the total error is extremely small. We have verified that boundary errors can be reduced by increasing the order and size of spacial domain and then machine precision can be reached again.
Investigating the effect of using approximate boundary conditions is shown in Fig. \ref{fig13}. Here we used the CW solution, (\ref{bseq}), to calculate the boundary conditions (\ref{req1b}) and (\ref{req2b})  with the same parameters as used for the dark soliton in Fig. \ref{fig9}. Similar to the bright soliton case, the error using approximate boundary conditions is almost indistinguishable from that with exact boundary conditions.

\section{comparison with other methods}
\label{compsec} Among the many numerical methods developed to solve
the NLSE, the so-called G-FDTD was shown to exceed  by orders of
magnitude the accuracy of all other methods
\cite{meth2,meth3,meth12,phd}. Therefore, we restrict the comparison
to this method. We start by pointing out the similarities and
differences between the G-FDTD and present method.

\begin{figure}[h]
    \centering
    \includegraphics[scale=0.62]{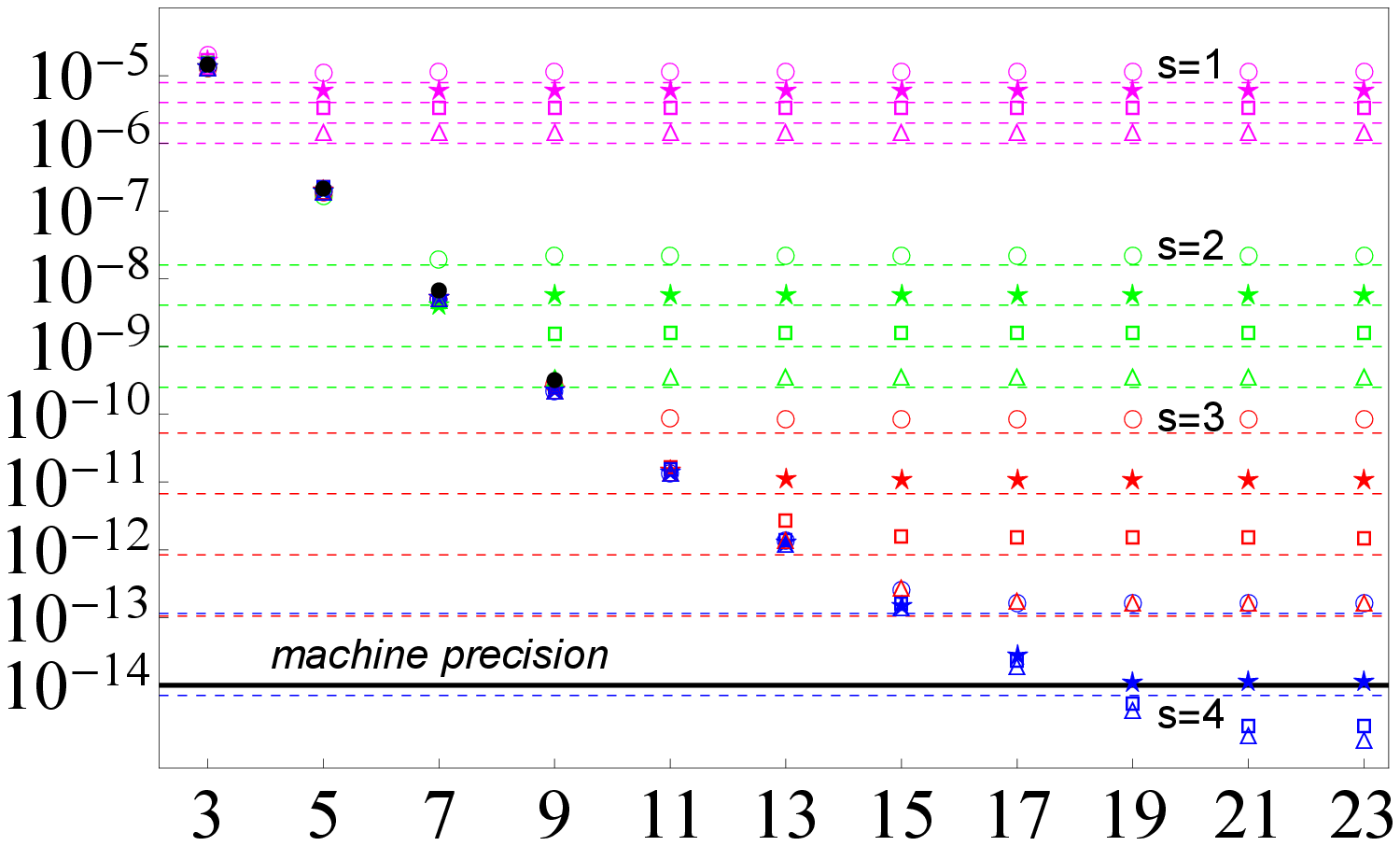}
\begin{picture}(5,5)(5,5)
\put(-277,70) {\rotatebox{90}{\scriptsize$\text{Log}(error)$}}
\end{picture}
\begin{picture}(5,5)(5,5)
\put(-127,-4) {\scriptsize$p$}
\end{picture}\vspace{0.2cm}
    \caption{ Error in the numerical solution of Eq. (\ref{nlse}) using the bright soliton (\ref{bs}) as an initial profile for different values of $p$, $s$, and $\Delta t$. Black filled circles correspond to $error_p$, as predicted by (\ref{errorp}) for $p=3, 5,7, 9$. Horizontal dashed lines correspond to $error_s$, as predicted by Eq. (\ref{errors}). The thick line near the bottom indicates the machine precision, $10^{-14}$. For each order the calculation is performed with $\Delta t$ (circles), $\Delta t/2$ (stars), $\Delta t/4$ (squares), and $\Delta t/8$ (triangles). Parameters: $\Delta t=0.001$, $n_x=500$, $L=40$, $\Delta x=L/(n_x-1)$, $k=4.0$, $x_0=0.0$, $A_0=1.0$, $g_1=-1$, $g_2=-2$, and for all points $t_f=n_t\times\Delta t=1$.}
    \label{fig3}
\end{figure}

\begin{figure}[h]
    \centering
    \includegraphics[scale=1.05]{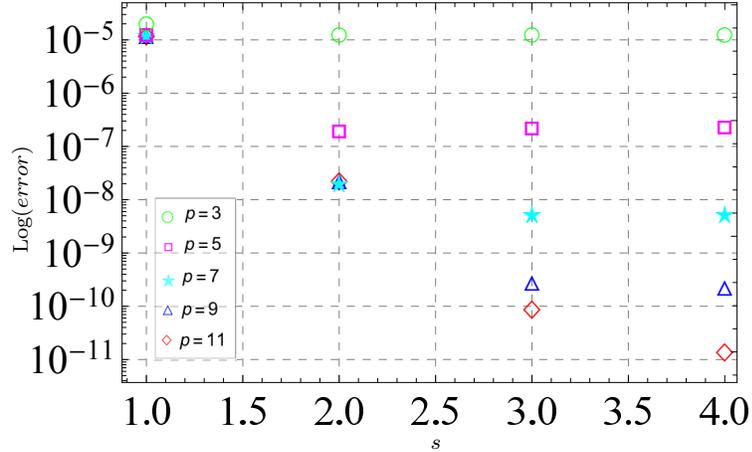}
    \begin{picture}(5,5)(5,5)
    \put(-277,70) {\rotatebox{90}{\scriptsize$\text{Log}(error)$}}
    \end{picture}
    \begin{picture}(5,5)(5,5)
    \put(-127,-4) {\scriptsize$s$}
    \end{picture}\vspace{0.2cm}
    \caption{ Error in the numerical solution of the bright soliton corresponding to Eq. (\ref{bright}) versus $s$ and different values of $p$ in semi-log scale.  Parameters: $\Delta t=0.001$, $n_x=500$, $n_t=1$, $L=40$, $\Delta x=L/(n_x-1)$, $k=4.0$, $x_0=0.0$, $A_0=1.0$, $g_1=-1$, $g_2=-2$.  The circles, squares, stars, triangles, and diamonds indicate the maximum error with $p=3,5,7,9,11$, respectively. For $p=13,15,17,19,21, 23$, the points, which are not shown here,  overlap with the diamonds $(p=11)$. }
    \label{fig4}
\end{figure}

\begin{figure}[h]
    \centering
    \includegraphics[scale=0.64]{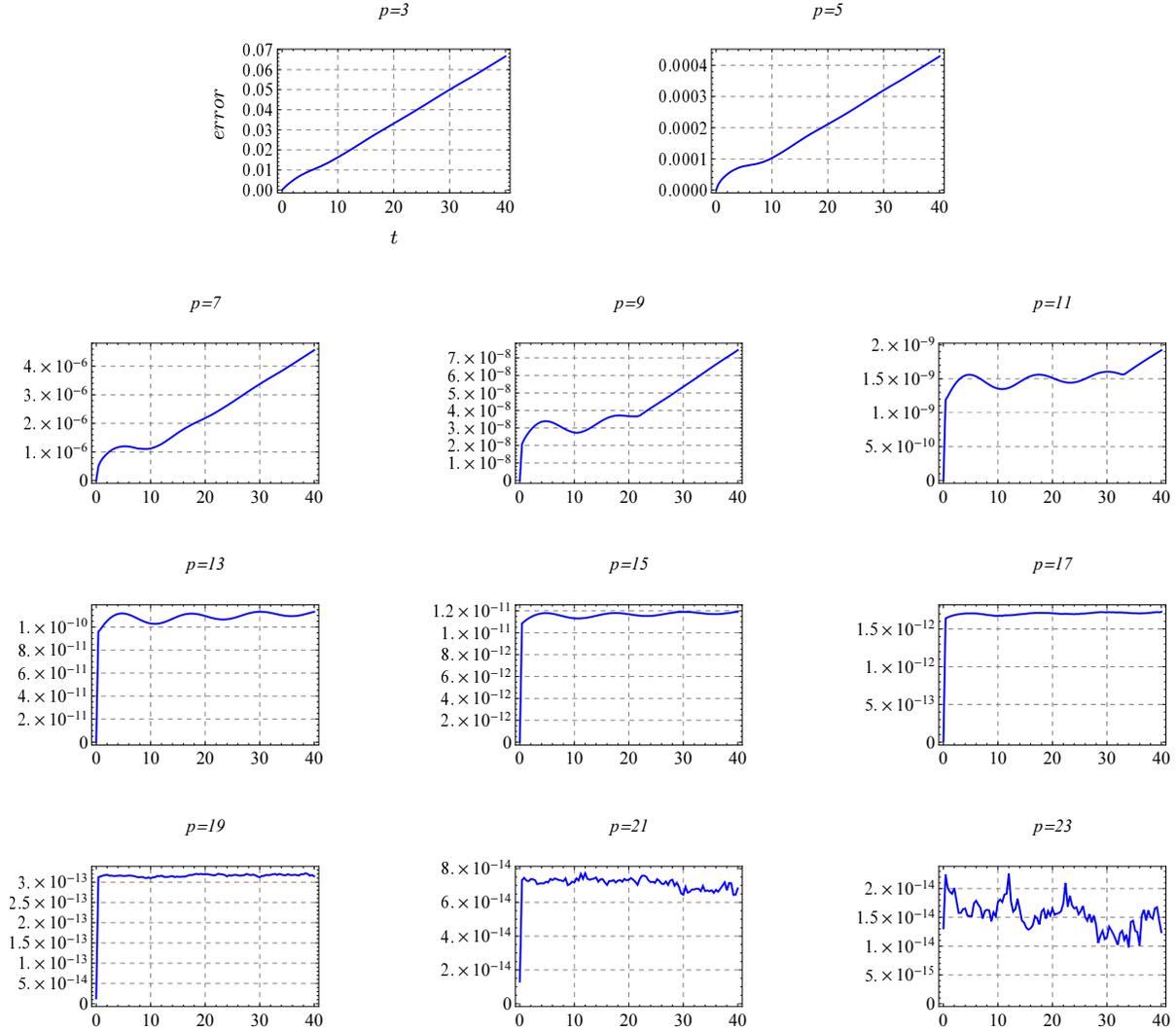}
    \begin{picture}(5,5)(5,5)
    \put(-370,350) {\rotatebox{90}{\scriptsize$error$}}
    \end{picture}
    \begin{picture}(5,5)(5,5)
    \put(-310,310) {\scriptsize$t$}
    \end{picture}
    \caption{Error evolution in the numerical solution of the bright soliton corresponding to (\ref{bright}).  Fourth order power series expansion is used for all cases. For $p\ge13$, the error saturates to a value independent of time. Machine precision is reached with $p=23$. Parameters: $n_x=8000$, $n_t=80000$, $L=800$, $\Delta t=0.0005$, $\Delta x=L/(n_x-1)$, $k=1.0$, $x_0=0.0$, $A_0=1.0$, $g_1=1/2$, $g_2=1$. }
    \label{fig5}
\end{figure}

\begin{figure}[h]
    \centering
    \includegraphics[scale=1.05]{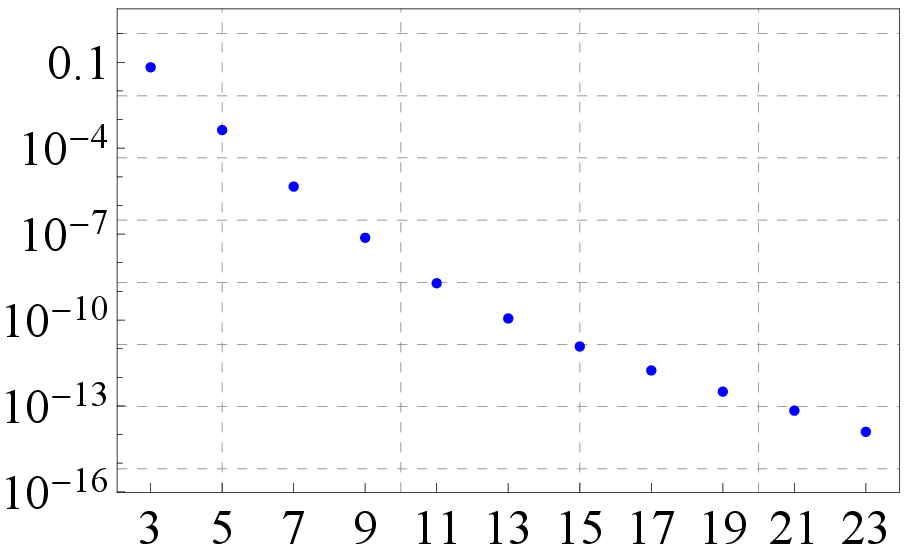}
    \begin{picture}(5,5)(5,5)
    \put(-283,67) {\rotatebox{90}{\scriptsize$\text{Log}(final\,\,error)$}}
    \end{picture}
    \begin{picture}(5,5)(5,5)
    \put(-128,-4) {\scriptsize$p$}
    \end{picture}\vspace{0.2cm}
    \caption{Error at the final evolution time for sub-figures of Fig. \ref{fig5} on a semi-log scale. For small values of $p$, the points approach the asymptote $\text{Log}(final\,\,error)=17.60-15.78\, \text{Log}(p)$, which gives: $final\,error=4.40\times10^7\,p^{-15.8}$.}
    \label{fig6}
\end{figure}

\begin{figure}[h]
    \centering
    \includegraphics[scale=1.08]{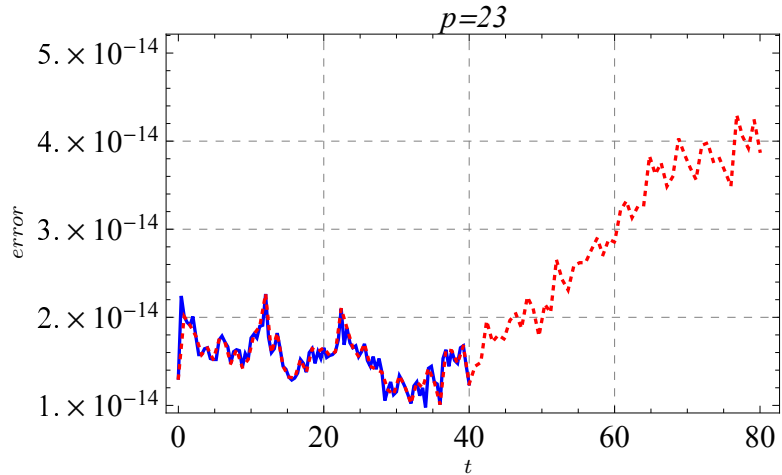}
    \begin{picture}(5,5)(5,5)
    \put(-289,77) {\rotatebox{90}{\scriptsize$error$}}
    \end{picture}
    \begin{picture}(5,5)(5,5)
    \put(-127,-4) {\scriptsize$t$}
    \end{picture}\vspace{0.2cm}
    \caption{Error evolution in the numerical solution of the bright soliton corresponding to Eq.(\ref{bright}). Dotted red curve corresponds to $p=23$ of Fig. \ref{fig5} with double evolution time.  Fourth order power series expansion is used in this case. The corresponding curve from Fig. \ref{fig5} is re-plotted here, with solid blue, and can be distinguished in the time interval $[0,40]$.   Parameters: $n_x=8000$, $n_t=160000$, $L=800$, $\Delta t=0.0005$, $\Delta x=L/(n_x-1)$, $k=1.0$, $x_0=0.0$, $A_0=1.0$, $g_1=1/2$, $g_2=1$.}
    \label{fig7}
\end{figure}

\begin{figure}[h]
    \centering
    \includegraphics[scale=0.64]{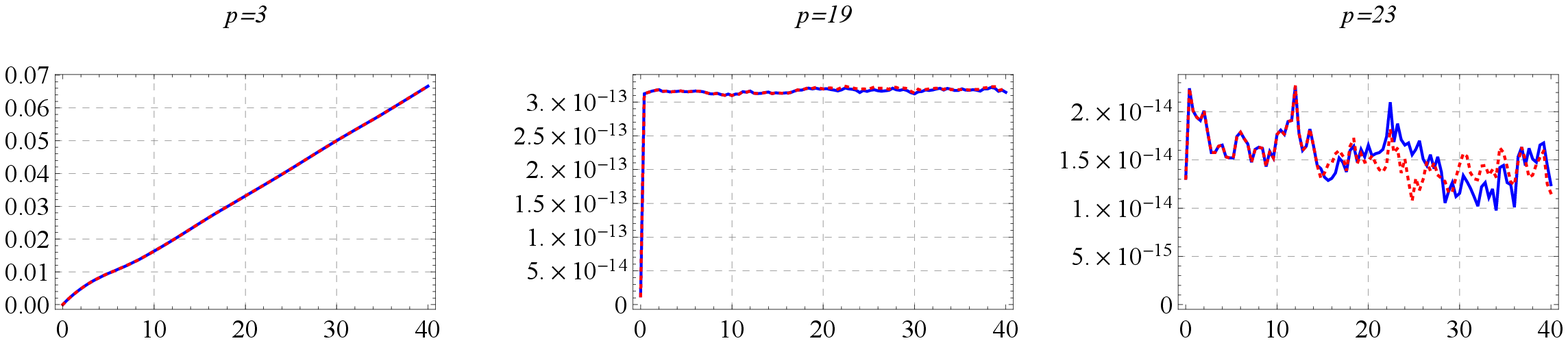}
    \begin{picture}(5,5)(5,5)
    \put(-430,38) {\rotatebox{90}{\scriptsize$error$}}
    \end{picture}
    \begin{picture}(5,5)(5,5)
    \put(-366,-5) {\scriptsize$t$}
    \end{picture}\vspace{0.2cm}
    \caption{Error evolution in the numerical solution of the bright soliton corresponding to Eq. (\ref{bright}). Dotted red curve corresponds to approximate boundary conditions. Solid blue curve corresponds to exact boundary conditions.  Fourth order power series expansion is used for all cases. Parameters used are the same of those in Fig. \ref{fig5}.}
    \label{fig8}
\end{figure}

\begin{figure}
    \centering
    \includegraphics[scale=0.64]{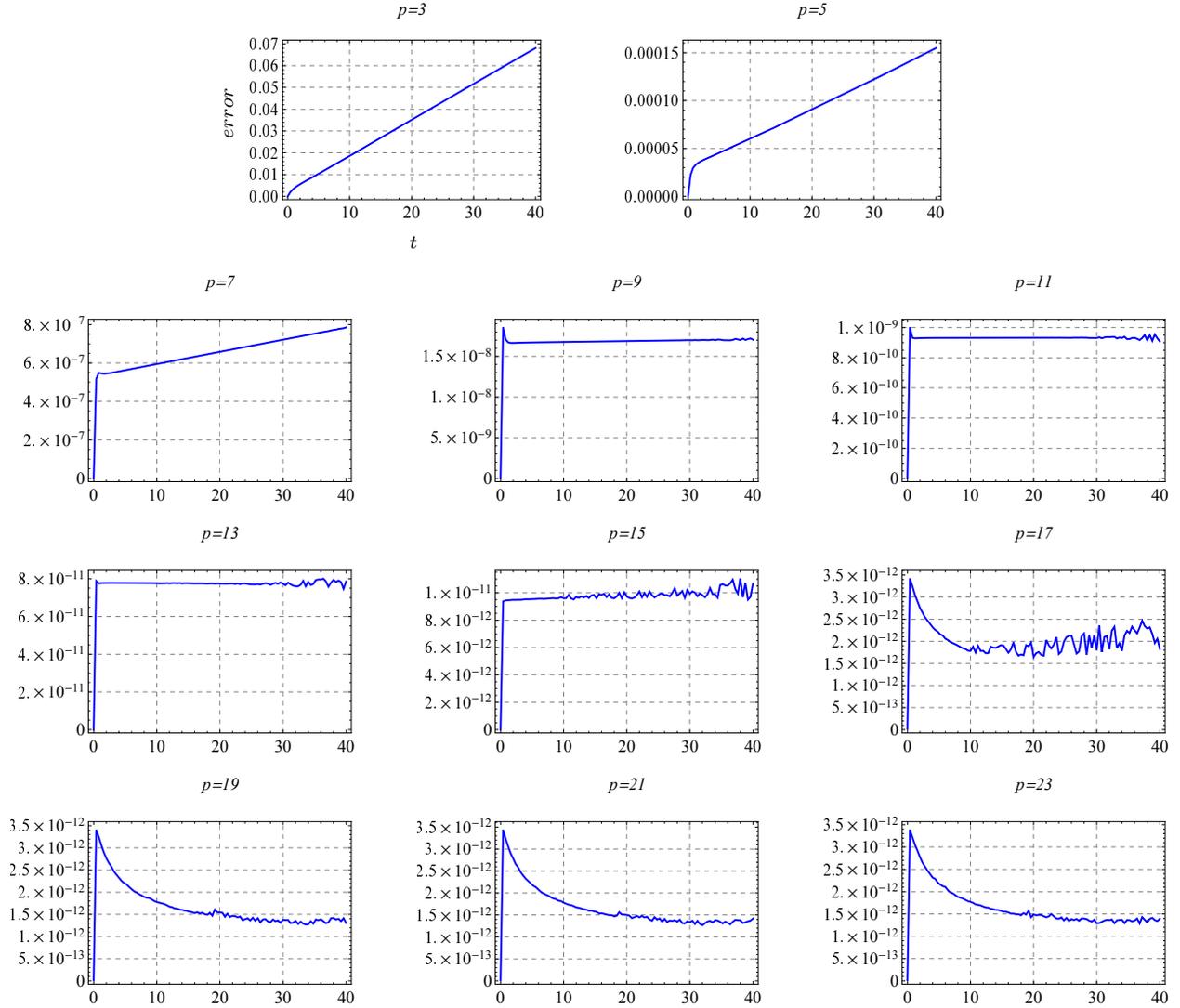}
    \begin{picture}(5,5)(5,5)
    \put(-135,370) {\rotatebox{90}{\scriptsize$error$}}
    \end{picture}
    \begin{picture}(5,5)(5,5)
    \put(-68,325) {\scriptsize$t$}
    \end{picture}
    \caption{Error evolution in the numerical solution of the dark soliton corresponding to Eq.(\ref{dark}) using increasing number of points in the second derivative discretization.  Fourth order power series expansion is used for all cases. For $p>7$, the error saturates at a value independent of $t$. The saturation value decreases with $p$.  Parameters: $n_x=4000$, $n_t=80000$, $L=400$, $\Delta t=0.0005$, $\Delta x=L/(n_x-1)$, $k=1.0$, $x_0=0.0$, $A_0=1.0$, $g_1=1/2$, $g_2=-1$. }
    \label{fig9}
 \end{figure}

\begin{figure}
    \centering
    \includegraphics[scale=1.05]{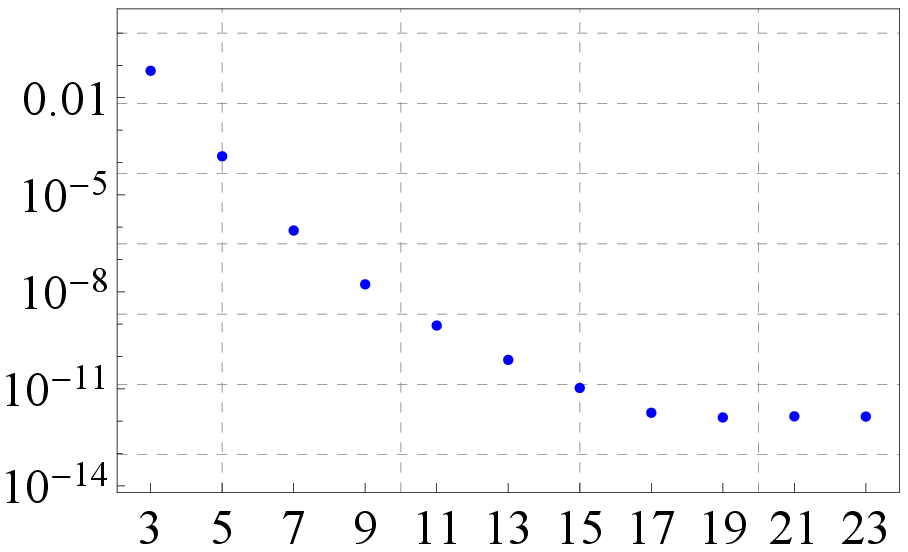}
        \begin{picture}(5,5)(5,5)
    \put(-282,60) {\rotatebox{90}{\scriptsize$\text{Log}(final\,\,error)$}}
    \end{picture}
    \begin{picture}(5,5)(5,5)
    \put(-128,-4) {\scriptsize$p$}
    \end{picture}\vspace{0.2cm}
    \caption{Error at the final evolution time for sub-figures of Fig. \ref{fig9} on a semi-log scale. For small values of $p$, the points approach the asymptote $\text{Log}(final\,\,error)=-6.22972-6.97164\, \text{Log}(p)$, leading to $final\,error=0.00197\,p^{-6.9716}$.}
    \label{fig10}
\end{figure}

\begin{figure}
    \centering
    \includegraphics[scale=1.05]{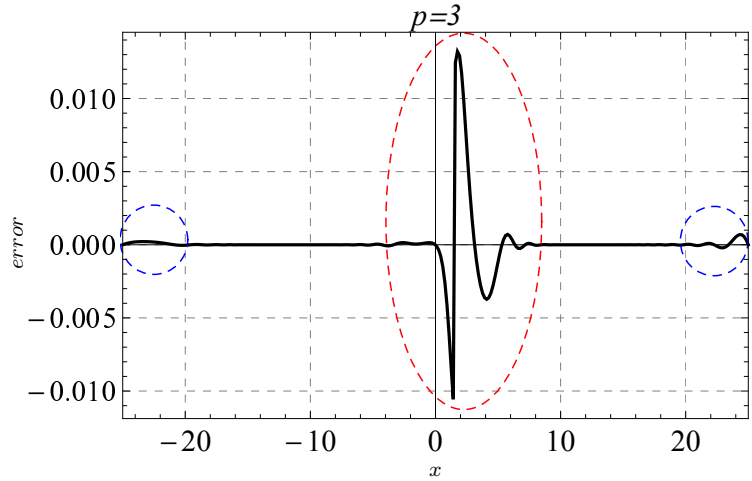}
    \begin{picture}(5,5)(5,5)
    \put(-278,74) {\rotatebox{90}{\scriptsize$error$}}
    \end{picture}
    \begin{picture}(5,5)(5,5)
    \put(-129,-4) {\scriptsize$x$}
    \end{picture}\vspace{0.2cm}
    \caption{Snapshot of error in the numerical solution of the dark soliton given by Eq. (\ref{dark}) at $t=1.5$. The circles in dashed blue show the source of  error from the edges and the circle in dashed red shows the source of error from the center. Parameters: $p=3$, $s=4$, $n_x=300$, $n_t=1000$, $L=50$, $\Delta t=0.01$, $\Delta x=L/(n_x-1)$, $k=1.0$, $x_0=0.0$, $A_0=1.0$, $g_1=1/2$, $g_2=-1$.}
    \label{fig11}
\end{figure}

\begin{figure}
    \centering
    \includegraphics[scale=1.00]{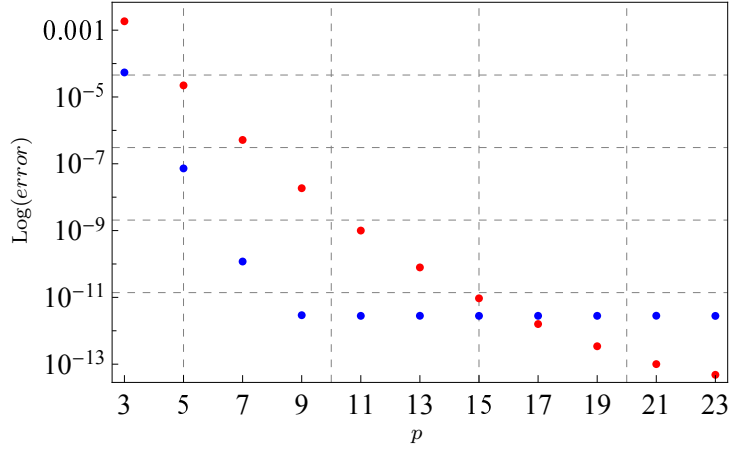}
    \begin{picture}(5,5)(5,5)
\put(-270,70) {\rotatebox{90}{\scriptsize$\text{Log}(error)$}}
\end{picture}
\begin{picture}(5,5)(5,5)
\put(-127,-4) {\scriptsize$p$}
\end{picture}\vspace{0.2cm}
\caption{Error at the edge and center of the spatial grid plotted in a semi-log scale. Red and blue points correspond to error at the centre and edges of the spacial grid, respectively. Each point corresponds to the error in a subfigure of Fig. \ref{fig9} for the same value of $p$ and $t=40$. Parameters used are the same of those in Fig. \ref{fig9}.}
\label{fig12}
\end{figure}

\begin{figure}
    \centering
    \includegraphics[scale=0.62]{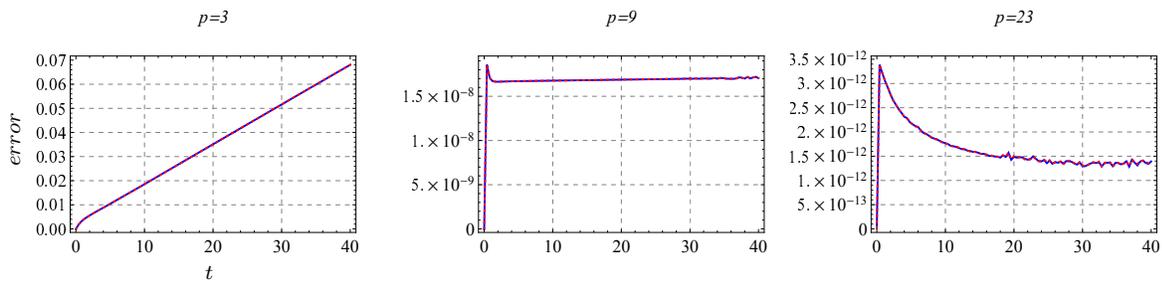}
    \begin{picture}(5,5)(5,5)
    \put(-433,38) {\rotatebox{90}{\scriptsize$error$}}
    \end{picture}
    \begin{picture}(5,5)(5,5)
    \put(-368,-5) {\scriptsize$t$}
    \end{picture}\vspace{0.2cm}
    \caption{Error evolution in the numerical solution of the dark soliton given by Eq.(\ref{dark}). The dotted red curve corresponds to approximate boundary conditions. The solid blue curve corresponds to exact boundary conditions.  Fourth order power series expansion is used for all cases. Parameters used are the same of those in Fig. \ref{fig9}.}
    \label{fig13}
\end{figure}

Both methods use a power series expansion to integrate the time operator. However, in the G-FDTD method, the Crank-Nicolson time stepping method is used and the expansion is restricted to odd powers of $\Delta t$. Here, we do not use the Crank-Nicolson time stepping  and we include all powers in the time expansion, as shown in (\ref{seru}) and (\ref{serv}). The use of the Crank-Nicolson method reduces the error in the time evolution by one order of magnitude in $\Delta t$. This additional accuracy comes on the expense of memory and run time cost; the evolution to $t+\Delta t$ requires the knowledge of the profiles at $t$ and $t+\Delta t/2$ and this has to be done for the real and imaginary parts of the profile. Thus, it requires at least 4 times run time and memory size compared with typical time stepping. This additional memory and computing time cost will increase dramatically in higher dimensions. On the other hand, our analysis of the error in Fig. \ref{fig3} has shown that, before the saturation region is reached, the error from the time evolution, $error_s$, is orders of magnitude smaller than the error from the spacial discretisation of the second derivative, $error_p$. Therefore, the additional accuracy brought by the use of the Crank-Nicolson time stepping is really not needed at this stage. It may have an advantage in case very long time evolution is needed, but even in that case, increasing the order, $s$, will lead to the required accuracy with less run time and memory storage.

In the G-FDTD method, an approximation was used in the calculation of the time evolution of the power expansion coefficients which is to consider the terms $|\psi(x,t)|$  as constants. This may not be justifiable with rapid and steep changes in $\psi(x,t)$. Here, we do not make this approximation. Our approach for the calculation of the coefficients of the power series is slightly different than that  of Ref. \cite{meth2,meth3,meth12,phd} allowing us to calculate the recursion relations and coefficients of the power series without any approximation.

Another point of similarity is that higher order discretisation of the second spacial derivative is used. In Ref. \cite{meth2,meth3,meth12,phd}, only the 3- and 5-point central difference point are used. Here, we use mostly up to the 23-point central difference formula to show that machine accuracy can indeed be reached with such a procedure. Our procedure allows for a straightforward and rather easy implementation of larger number of points.

The treatment of boundary conditions has also similarities and differences. Both methods use the exact localised initial solution or the CW solution to calculate the time evolution of the boundary points. In Ref. \cite{meth3,phd}, the first and last 6 points ($p=5$ is used) were set to have the value of the exact solution for all times. In the present method, we set only the first and last $2$ points to the exact values. The difference in number of boundary points is due to the difference in procedure of calculating the coefficients of the time power series.  We believe it is more accurate to set only $(p-1)/2$ boundary points and not $p+1$, since the $p$-central difference formula correlates only the left or right $(p-1)/2$ points to the central point. The effect of this difference in handling the boundary conditions will have an evident effect on solutions with finite background, such as dark and peregrine solitons, as we will point out below. \\\\

\textbf{\textit{Bright soliton:}}\\
In the following, we consider the same parameters taken by Ref.
\cite{meth3,phd} to calculate the accuracy versus the exact solutions, namely  $g_1=-1,\,g_2=-2, A_0=1,\, k=4,\,x_0=-10$. Since in Ref. \cite{meth3,phd}, the error is calculated as $error=\sqrt{(1/n_x)\sum_{i=1}^{n_x}\left(|\psi_{\rm numerical}(x_i,t_f)|-|\psi_{\rm exact}(x_i,t_f)|\right)^2}$, we use this definition in this and the dark soliton comparison. Table \ref{table1} shows that similar errors are obtained for $n_x=100,\,200,\,300$, but not for $n_x=400$ where our method gives about $40\%$ smaller error. To understand this difference in error, we calculate the convergence rate defined by $R=\log(error_2/error_1)/\log(n_{x2}/n_{x1})$, where $error_{1,2}=error(n_{x1,2})$ calculated at two values of $n_x$. Convergence rate gives the exponent in the power law dependence of error on $\Delta x$, namely $error\propto (\Delta x)^{R}$. For the fourth order central point formula, used here, the error is proportional to $(\Delta x)^4$ and hence the convergence rate should be $R=4$. Table \ref{table1} shows that indeed for both methods $R\approx4$ for $n_x=100,\,200,\,300$, but it is not so for the G-FTDT method with $n_x=400$. To reach the theoretical convergence rate of 4, the value of $\Delta t$ needs to be reduced which requires increasing $n_t$ by the same factor so that the final time remains unchanged. This will of course require increased run time by the same factor.  In Table \ref{table2}, we recalculate the error using $\Delta t=10^{-6}$.  Error and convergence rate of the present method are almost the same as those with $\Delta t=10^{-4}$, which is understood since the error at this stage is dominated by $error_p$. On the other hand, the error of the G-FDTD method at $n_x=400$ has decreased and is now indeed close to that of the present method with a convergence rate approaching 4. The CPU time in this case is $170.7$ s which is to be compared with $0.6$ s for the present method at $n_x=400$ in Table
\ref{table1}. Thus, our code is more than 100 times faster than that of the G-FDTD method for about the same error and convergence rate.  The difference in speed grows with larger $n_x$. Data for the G-FDTD code in Table \ref{table1} is not shown for $n_x>400$, since it is unstable with $\Delta t=10^{-4}$, while our code continues to be stable for much larger values of $n_x$ with increasing accuracy and keeping the convergence rate approaching 4. With $\Delta t=10^{-6}$, the G-FTDT code is stable for $n_x>400$, but with convergence rate deviating from the theoretical value of 4 for larger $n_x$. Again, this can be fixed by decreasing $\Delta t$ further which will require more run time. It should be noted that we do not use in the G-FDTD code the additional refinement, used in Ref. \cite{phd}, of reducing the error to the machine precision from one step to the other, in order to preserve the norm and energy. This would significantly slower the G-FTDT code even further.

To show the high potential of the present method, we repeat the above calculations with larger number of central point formula. In Table \ref{table3}, we show the results for $p=11,\,15,\,23$. The significant reduction in error is obvious with CPU times on the order of 1 s. The error can be seen to drop easily down to machine precision with $p=23$ and $n_x>350$. The convergence rates for $p=11$ and $p=15$ are close to the theoretical values of 10 and 14, respectively. However, for $p=23$, the convergence rate reaches a maximum of 19.2 at $n_x=400$ and then starts to drop. Unlike the similar case above with the G-FDTD method, it will not be possible to increase the convergence rate to the theoretical value of 22 by decreasing $\Delta t$ since the error at this stage has reached the machine precision and does not any more depend on the parameters of the method. Remarkably, machine accuracy is reached with a CPU time less than 2 s. 

A high accuracy method based on wavelets expansion, was developed in Ref. \cite{meth7}. The present example compares with Example 2 in that reference. For $\Delta t=0.01$ and $n_x=200$, an error of $8.94\times10^{-5}$ at $t=1.0$ was obtained with a convergence rate close to 7. In the present method, this matches $s=8$. Since we do not consider this value, we compare with $s=7$ and $s=9$, which have theoretical convergence rates 6 and 8, respectively. The error in these two cases turn out to be $1.05\times10^{-4}$ and $1.16\times10^{-5}$, respectively. The CPU times are 0.01 s and 0.5 s for the former and latter cases, respectively. If we take the average of two error values for $s=7$ and $s=9$, the error value of the present method  will be about 30\% less  than that of the Ref. \cite{meth7}. Unfortunately, comparison of CPU time is not possible because they it is not reported in that reference. 

\begin{figure}[h]
    \centering
    \includegraphics[width=15cm]{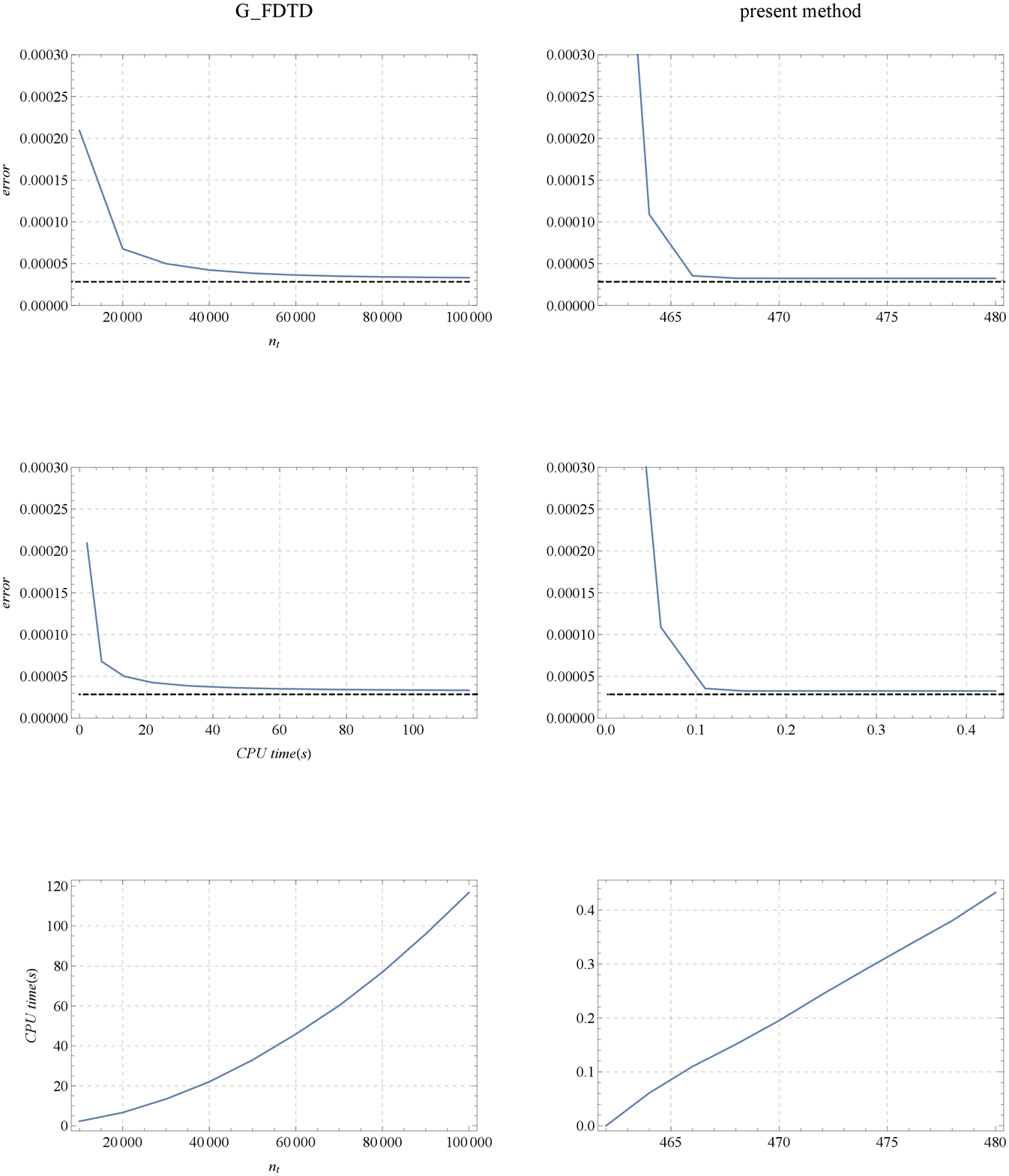}
    \caption{Error and CPU time of the numerical bright soliton solution  obtained by the present and G-FDTD methods. Error is calculated using the exact analytical bright soliton solution (\ref{bs}). Horizontal dashed asymptote corresponds the theoretical estimate given by (\ref{errorp2}). Parameters: $g_1=-1,\,g_2=-2, A_0=1,\, $t$,\,x_0=-10$, $p=5,\,s=3$, $n_x=500$. }
    \label{fig14}
\end{figure}

Another important feature to  present  is the saturation of error at a constant minimum independent of $\Delta t$. This can be anticipated in view of our discussion of Fig. \ref{fig3}. Error will generally reduce with decreasing $\Delta t$, but when $\Delta t$ is small enough, the error from time stepping will be less than that of the central point formula and thus the total error will be independent of $\Delta t$. The error in this case will be determined essentially by Eq. (\ref{errorp}) but without the operator $\Delta t\,\partial/\partial t$ since the error at this stage is dominated by that of the central point formula
\begin{eqnarray}
error_p&=&\max\left[\frac{(\Delta x)^{p-1}}{(p+1)!}\frac{\partial |\psi^{p+1}(x,t)|}{\partial x^{p+1}}\right].\label{errorp2}
\end{eqnarray}
It is important to know the maximum $\Delta t$ in this region for which the error is independent of $\Delta t$  in order to save on the run time; no need to run the code with a very small, and hence time consuming $\Delta t$, while a larger value can produce the same error with less run time. In Fig. \ref{fig14}, this is manifested through a number of plots. At first, the upper panel shows that error of the present method reduces rather sharply to its saturation level, as predicted by (\ref{errorp2}), around $n_t\sim465$ ($\Delta t=1/n_t\approx2.15\times10^{-3}$) over a range of $n_t=2$. On the other hand, the error in the  G-FDTD method reduces in a slower rate to the same saturation level  at around $n_t=20000$ ($\Delta t=1/n_t\approx5\times10^{-5}$). The middle panel shows the huge difference in CPU times between the two methods. The bottom panel shows that while CPU time of the present method grows linearly with $n_t$, it grows quadratically with $n_t$ for the G-FDTD method.\\\\

\textbf{\textit{Dark soliton:}}\\
In Table \ref{table4}, we show that the accuracy of the present method keeps increasing with increasing $n_x$ while that of the G-FDTD saturates at a certain value. As explained above, the convergence rate is still far from the theoretical value of 4 and thus smaller $\Delta t$ is needed in order to get smaller errors and better convergence rate.

\section{Other examples}
\label{othersec}
Here we put the method under tests of cases involving high curvatures and fast time evolution, namely the Peregrine  soliton and the two soliton molecule. The Peregrine soliton is characterised by high curvature at the time of its maximum peak. This will test the accuracy in the $p$-points formula for the second derivative. It will also test our treatment of the boundary conditions since the background for this soliton is finite. The two-soliton molecule is characterised by fast dynamics in the case when the two solitons coalesce. This will test the accuracy in the time power series method that integrates the time operator. Another important feature in both of these two solutions is that, unlike the previous two examples, their time evolution is nontrivial. In moving bright and dark solitons, the internal structure does not change. For  the current examples, the internal structure changes with time, which will result in larger errors, as we will see below.\\\\

\begin{figure}
    \centering
    \includegraphics[scale=1.08]{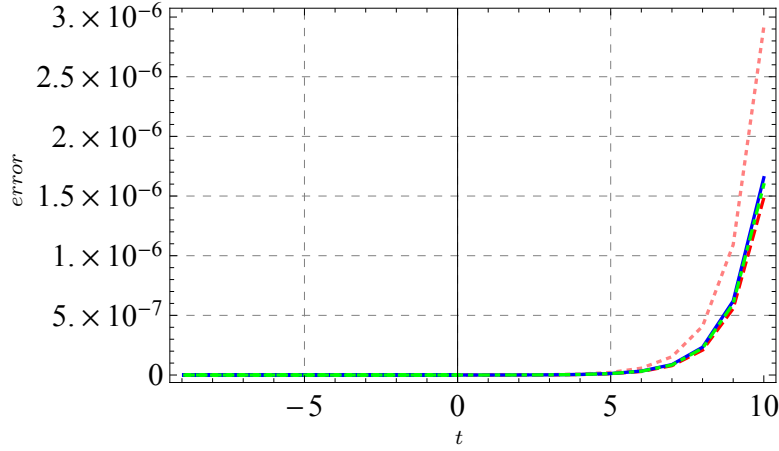}
   \begin{picture}(5,5)(5,5)
   \put(-288,85) {\rotatebox{90}{\scriptsize$error$}}
    \end{picture}
    \begin{picture}(5,5)(5,5)
    \put(-129,-5) {\scriptsize$t$}
    \end{picture}\vspace{0.2cm}
    \caption{Error evolution in the numerical solution of the Peregrine soliton given by Eq. (\ref{per}) using increasing number of points in the second derivative discretization. Dotted pink curve corresponds to $p=9$, dashed red curve corresponds to $p=23$, solid black curve corresponds to $p=25$, and green dotted-dashed green curve corresponds to $p=27$.  Fourth order power series expansion is used for all cases.  Parameters: $n_x=2000$, $n_t=200000$, $L=50$, $\Delta t=0.0001$, $\Delta x=L/(n_x-1)$, $g_1=1/2$, $g_2=1$, $t_0=-10$.}
    \label{fig15}
\end{figure}

\textbf{\textit{Peregrine soliton:}}\\
The exact Peregrine soliton of  Eq. (\ref{nlse}) takes the following expression
\begin{equation}\label{per}
\psi(x,t)=\frac{1}{\sqrt{g_2}}\,\Big[\frac{4+i\,8\,(t-t_0)}{1+4\,(t-t_0)^2+\frac{2}{{g_1}}\,(x-x_0)^2}-1\Big]\,e^{i\,[t-t_0+\phi_0]},
\end{equation}
where $g_2>0$.
The initial profile is started at $t=-10$ and evolved till $t=10$.
Figure \ref{fig15} shows the maximum error versus time for four values of $p$. The figure shows that the Peregrine soliton is highly demanding computationally, as hinted above. The error curves of $p=23,\,25,\,27$ are almost the same. Increasing $p$ will thus not enhance on the accuracy. Similar to the dark soliton, this terminal error is due to the finite background. It can be reduced by increasing the size of the spacial grid and increasing the order, $s$.\\\\

\begin{figure}
    \centering
    \includegraphics[scale=0.55]{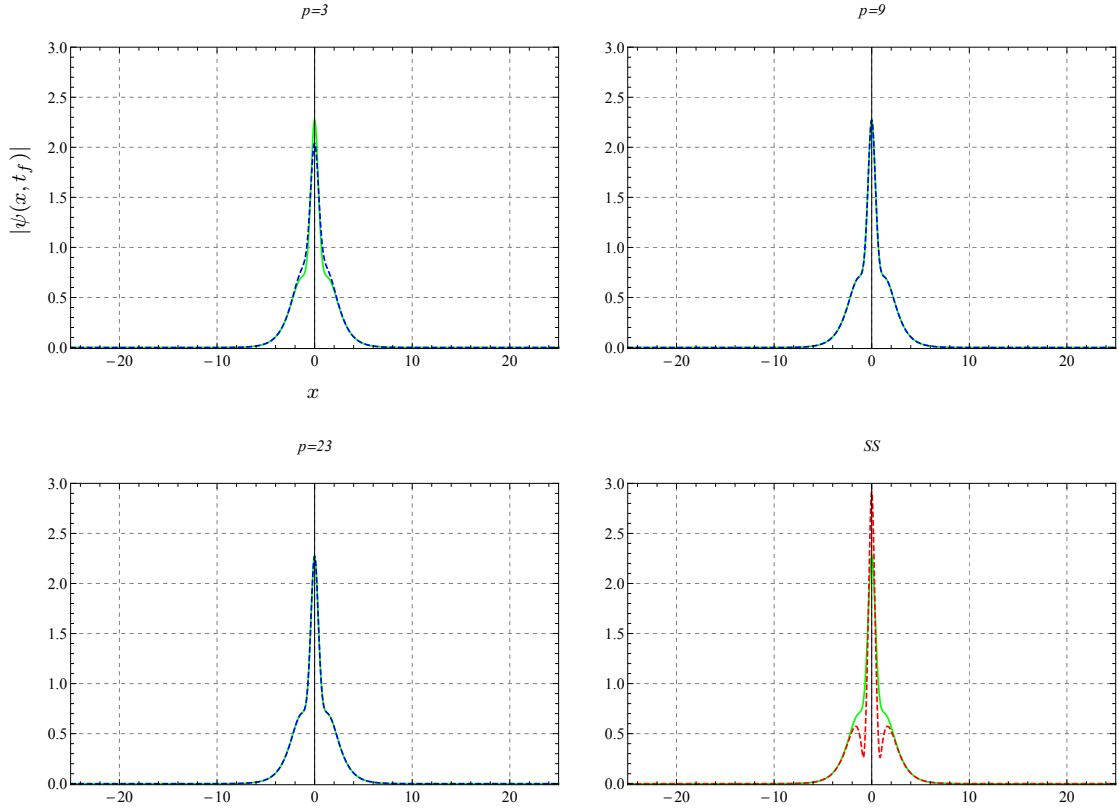}
    \begin{picture}(5,5)(5,5)
    \put(-417,220) {\rotatebox{90}{\scriptsize$|\psi(x,t_f)|$}}
    \end{picture}
    \begin{picture}(5,5)(5,5)
    \put(-313,158) {\scriptsize$x$}
    \end{picture}
    \caption{Profile of the  two-bright soliton solution profile at $t=80$. Solid green curve is the analytical exact solution given by Eq.(\ref{two}), dashed blue curves correspond to the present method using $p=3,9,23$, and dashed red curve corresponds to  the SS method.  Fourth order power series expansion is used for all cases with $n_x=4000$.  The dashed red curve is the result of the SS code with $n_x=4096$. Other parameters: $n_t=800000$, $L=50$, $\Delta t=0.0001$, $t_f=n_t\times\Delta t$, $\Delta x=L/(n_x-1)$, $g_1=1/2$, $g_2=1$, $t_0=\phi_{01}=\phi_{02}=\nu_1=\nu_2=x_{01}=x_{02}=0$, $\alpha_1=1$, $\alpha_2=2$.}
    \label{fig16}
\end{figure}

\begin{figure}
    \centering
    \includegraphics[scale=0.55]{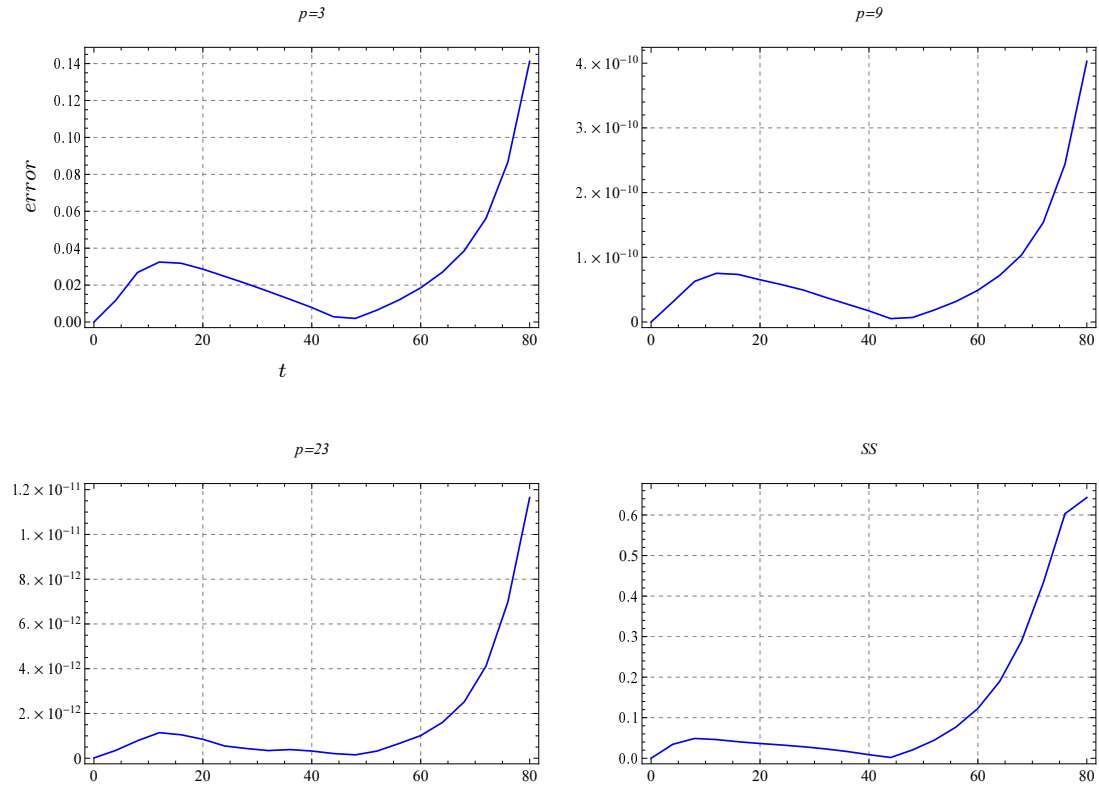}
    \begin{picture}(5,5)(5,5)
    \put(-403,220) {\rotatebox{90}{\scriptsize$error$}}
    \end{picture}
    \begin{picture}(5,5)(5,5)
    \put(-316,158) {\scriptsize$t$}
    \end{picture}
    \caption{Error evolution in the numerical solution of the two-bright soliton solution given by Eq. (\ref{two}) using $p=3,9,23$ in the second derivative discretization.  Fourth order power series expansion is used for all cases.  The last subfigure is the result of the SS method. Parameters used are the same as those in Fig. \ref{fig16}.}
    \label{fig17}
\end{figure}

\textbf{\textit{Two-bright soliton:}}\\
The two-bright soliton of  Eq. (\ref{nlse}) takes the form
\begin{equation}\label{two}
\psi(x,t)=\frac{1}{\sqrt{g_2}}\,[\psi_1(x,t)+\psi_2(x,t)],
\end{equation}
where \\\\
$\psi_1(x,t)=\frac{M_{12}\,[\gamma_1^{-1}(x,t)+\gamma_2^*(x,t)]-M_{22}\,[\gamma_2^{-1}(x,t)+\gamma_2^*(x,t)]}{M_{12}\,M_{21}\,[\gamma_1^*(x,t)+\gamma_2^{-1}(x,t)]\,[\gamma_1^{-1}(x,t)+\gamma_2^{*}(x,t)]-M_{11}\,M_{22}\,[\gamma_1^{-1}(x,t)+\gamma_1^{*}(x,t)]\,[\gamma_2^{-1}(x,t)+\gamma_2^{*}(x,t)]}$,\\\\
$\psi_2(x,t)=\frac{-M_{11}\,[\gamma_1^{-1}(x,t)+\gamma_1^*(x,t)]+M_{21}\,[\gamma_1^{*}(x,t)+\gamma_2^{-1}(x,t)]}{M_{12}\,M_{21}\,[\gamma_1^*(x,t)+\gamma_2^{-1}(x,t)]\,[\gamma_1^{-1}(x,t)+\gamma_2^{*}(x,t)]-M_{11}\,M_{22}\,[\gamma_1^{-1}(x,t)+\gamma_1^{*}(x,t)]\,[\gamma_2^{-1}(x,t)+\gamma_2^{*}(x,t)]}$,\\\\
$a_1>0$, $a_2>0$, $M_{jk}=
1/(\lambda_j+\lambda_k^*)$,
$\gamma_j(x,t)=e^{\frac{\lambda_j}{\sqrt{2\,g_1}}(x-x_{0j})+i\,[\lambda_j^2\,(t-t_0)/2+\phi_{0j}]}$, $\lambda_j=\alpha_j+i\,\nu_j$,\\
$\alpha_j$, $\nu_j$, $x_{0j}$, $t_0$, and $\phi_{0j}$ are arbitrary real constants.
Here, we compare our method with the split-step (SS) method and plot the numerical profiles of both methods together with the exact one in Fig. \ref{fig16}. While the profile of the present method is indistinguishable from the exact one, the profile of the SS method deviates significantly at large evolution times. It should be noted that we used $p=3,\,9,\,23$ for our method. The $p=3$ curve shows a slight deviation, but the $p=9$ and $p=23$, are almost identical to the exact profile. The error for these three values of $p$ and the SS code are plotted in Fig. \ref{fig17}. Clearly, the $p=9$ and $p=23$ cases give extreemly small errors for a long evolution time.

\section{Inhomogeneous NLSE}
\label{vsec} Here, we present a generalisation of the method to the
NLSE with an external potential. Then we consider an example of
soliton scattering by a reflections potential well with a soliton
speed close to the critical value for quantum reflection. The
outcome, in this case reflection or transmission, is very sensitive
to the accuracy of the numerical method used. The high accuracy
provided by the present method is crucial for obtaining the correct
scattering outcome and accounting accurately for the value of the
critical speed.

\begin{figure}
    \centering
    \includegraphics[scale=0.55]{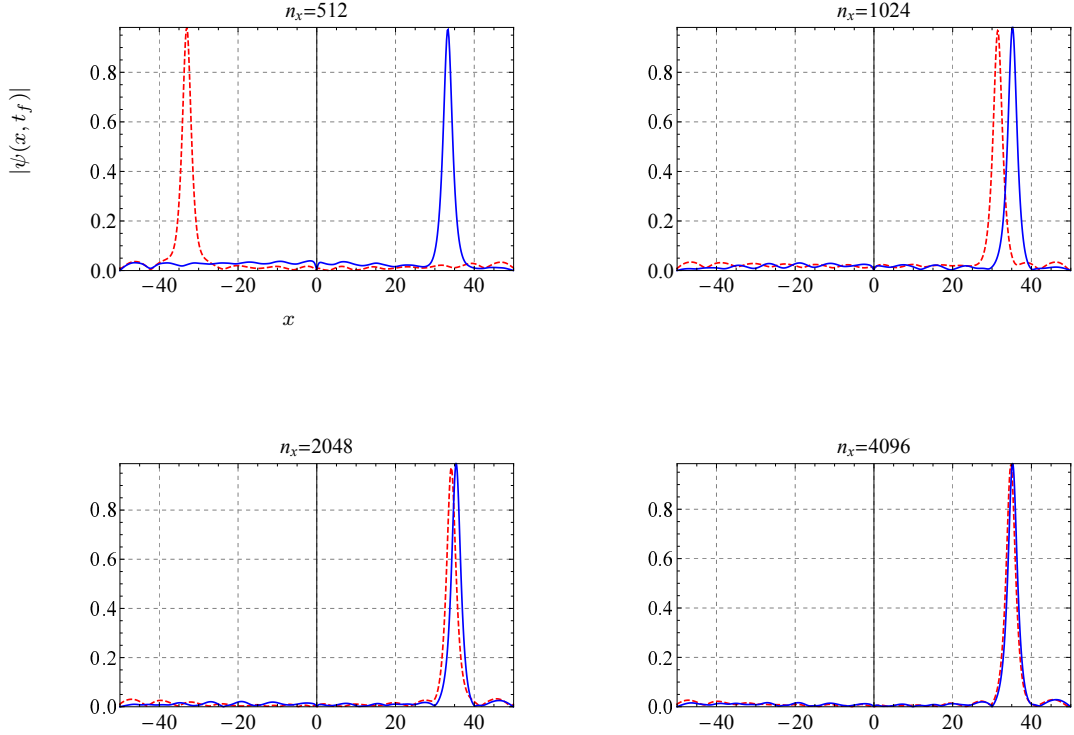}
    \begin{picture}(5,5)(5,5)
    \put(-400,215) {\rotatebox{90}{\scriptsize$|\psi(x,t_f)|$}}
    \end{picture}
    \begin{picture}(5,5)(5,5)
    \put(-305,158) {\scriptsize$x$}
    \end{picture}
    \caption{Snapshot of soliton profile after scattering by the reflectionless potential (\ref{pot}) at $t=153.186$. Solid blue curve corresponds to the present method and dashed red curve corresponds to the SS method. Other parameters: $\Delta t=0.1\times \Delta x^2$, $\Delta x=L/(n_x-1)$, $g_1=1/2$, $g_2=1$, $x_{0}=-10$, $k=0.331$, $\alpha=2$, $L=100$, $p=3$.}
    \label{fig18}
\end{figure}

\begin{figure}
    \centering
    \includegraphics[scale=0.22]{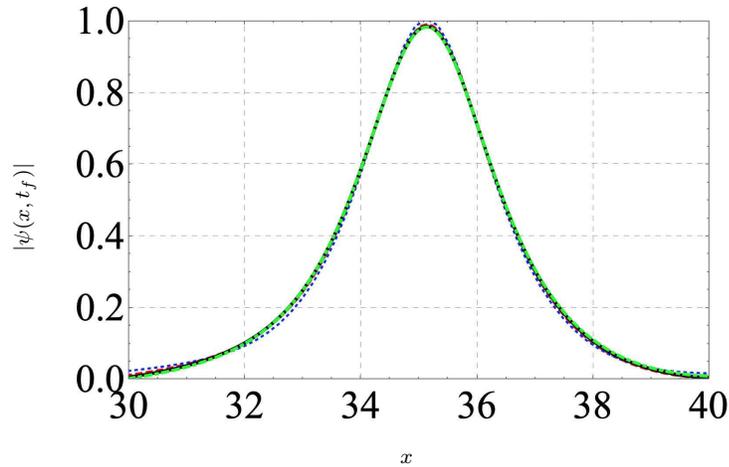}
     \begin{picture}(5,5)(5,5)
    \put(-269,77) {\rotatebox{90}{\scriptsize$|\psi(x,t_f)|$}}
    \end{picture}
    \begin{picture}(5,5)(5,5)
    \put(-131,-3) {\scriptsize$x$}
    \end{picture}\vspace{0.2cm}
    \caption{Snapshot of soliton profiles obtained by the present method after scattering by the reflectionless potential (\ref{pot}) at $t=153.186$. Dotted blue, dashed red, solid black, and dotted-dashed green curves correspond to $n_x=1/2\times1024,\,1024,\, 2\times1024, \, 4\times 1024$, respectively. Other parameters: $\Delta t=0.1\times \Delta x^2$, $\Delta x=L/(n_x-1)$, $g_1=1/2$, $g_2=1$, $x_{0}=-10$, $k=0.331$, $\alpha=2$, $p=23$.}
    \label{fig19}
\end{figure}

In the presence  of an external potential, $V(x)$, the NLSE can be expressed as
\begin{eqnarray}\label{V0}
i\, \frac{\partial}{\partial t}\psi(x,t)+g_1\,\frac{\partial^2}{\partial x^2}\psi(x,t)+g_2\,|\psi(x,t)|^2\,\psi(x,t)-V(x)\,\psi(x,t)&=&0.
\end{eqnarray}
Writing the general solution in the cartesian complex form $\psi(x,t)=u(x,t)+i\,v(x,t)$,  where $u(x,t)$ and $v(x,t)$ being real functions, and inserting in (\ref{V0}), generates two equations from the real and imaginary parts
\begin{eqnarray}
g_2\,u^3(x,t)+g_1\,\frac{\partial^2}{\partial x^2}u(x,t)+g_2\,u(x,t)\,v^2(x,t)-\frac{\partial}{\partial t}v(x,t)-V(x)\,u(x,t)&=&0,\\\label{V1}
g_2\,v^3(x,t)+g_1\,\frac{\partial^2}{\partial x^2}v(x,t)+g_2\,v(x,t)\,u^2(x,t)-\frac{\partial}{\partial t}u(x,t)-V(x)\,v(x,t)&=&0.\label{V2}
\end{eqnarray}
Substituting the power series expansions (\ref{uvexp1}) and (\ref{uvexp5}), the recursion relations are obtained as
\begin{eqnarray}
a_1&=&-g_2\,a_0^2\,b_0-g_2\,b_0^3+b_0\,V(x)-g_1\,b_0^{\prime\prime},\\
a_2&=&\frac{1}{2}\left[-2\,g_2\,a_0\,a_1\,b_0-g_2\,a_0^2\,b_1-3\,g_2\,b_0^2\,b_1+b_1\,V(x)-g_1\,b_1^{\prime\prime}\right],\\
a_3&=&\frac{1}{3}\Big[-g_2\,a_1^2\,b_0-2\,g_2\,a_0\,a_2\,b_0-2\,g_2\,a_0\,a_1\,b_1-3\,g_2\,b_0\,b_1^2\nonumber\\&&-g_2\,a_0^2\,b_2-3\,g_2\,b_0^2\,b_2+b_2\,V(x)-g_1\,b_2^{\prime\prime}\Big],\\
a_4&=&\frac{1}{4}\Big[-g_2\,a_1^2\,b_1-g_2\,b_1^3-2\,g_2\,a_0\,a_3\,b_0-2\,g_2\,a_0\,a_2\,b_1-6\,g_2\,b_0\,b_1\,b_2\nonumber\\&&-2\,g_2\,a_1\,a_2\,b_0-2\,g_1\,a_1\,a_0\,b_2-g_2\,a_0^2\,b_3-3\,g_2\,b_0^2\,b_3+b_3\,V(x)-g_1\,b_3^{\prime\prime}\Big],\\
b_1&=&g_2\,a_0\,b_0^2+g_2\,a_0^3-a_0\,V(x)+g_1\,a_0^{\prime\prime},\\
b_2&=&\frac{1}{2}\left[2\,g_2\,a_0\,b_0\,b_1+g_2\,b_0^2\,a_1+3\,g_2\,a_0^2\,a_1-a_1\,V(x)+g_1\,a_1^{\prime\prime}\right],\\
b_3&=&\frac{1}{3}\Big[g_2\,a_0\,b_1^2+2\,g_2\,a_0\,b_2\,b_0+2\,g_2\,a_1\,b_0\,b_1+3\,g_2\,a_0\,a_1^2\nonumber\\&&+g_2\,a_2\,b_0^2+3\,g_2\,a_0^2\,a_2-a_2\,V(x)+g_1\,a_2^{\prime\prime}\Big],\\
b_4&=&\frac{1}{4}\Big[g_2\,a_1\,b_1^2+g_2\,a_1^3+2\,g_2\,a_0\,b_3\,b_0+2\,g_2\,a_0\,b_1\,b_2+6\,g_2\,a_0\,a_1\,a_2\nonumber\\&&+2\,g_2\,a_1\,b_2\,b_0+2\,g_2\,a_2\,b_0\,b_1+g_2\,b_0^2\,a_3+3\,g_2\,a_0^2\,a_3-a_3\,V(x)+g_1\,a_3^{\prime\prime}\Big].
\end{eqnarray}
The boundary conditions are treated here in a similar manner as in the homogeneous case, namely using Eqs. (\ref{req1b}) and (\ref{req2b}). Finally, the time evolution is determined by the algorithm of the previous section but using the above modified recursion relations.\\\\
\textbf{\textit{Example: Soliton scattering by a reflectionless potential well}}\\
The bright soliton described by (\ref{bright}) is scattered by the following reflectionless potential well
\begin{eqnarray}\label{pot}
V(x)=-\frac{V_0^2}{\text{cosh}^2(\alpha\,x)},
\end{eqnarray}
where $V_0$ and $\alpha$ being arbitrary real constants. It is established that below a critical speed, the soliton will reflect. This is known as quantum reflection since it occurs due to a repulsive force of  interaction between a trapped mode formed from the tail of the incoming soliton with the rest of the soliton. We use a soliton speed very close to the critical value and observe the outcome in terms of accuracy of the method. We also compare our results with the SS method. In Fig. \ref{fig18}, we plot the soliton profiles long after scattering by the potential. For the crudest version of our code, namely $p=3$ and $n_x=512$, the soliton transmits. Using the same parameters, the SS code leads to reflection. Since an exact analytical solution is not available, we use the comparative analysis to have an estimate on the accuracy of our results. We increase $n_x$ till the profile saturates at a certain shape. Considering four values of $n_x$, Fig. \ref{fig19} shows that the profile of our code is saturating in the transmission region at around $x=35$. The profile of the SS code transfers from reflection to the transmission region and gradually approaches the profile of the present method. This is shown more clearly where we use $p=23$ to see that saturation is already reached where the shape and position of the profiles are the same. There are small deviations for the $n_x=512$ case but they gradually disappear with increasing $n_x$.

\section{Two-coupled NLSE}
\label{twosec} Here we apply the method to two-coupled NLSE and use
the dark-bright soliton exact solution to check the accuracy.  The
two-coupled NLSE reads
\begin{eqnarray}
i\, \frac{\partial}{\partial t}\psi_1(x,t)+g_{10}\,\frac{\partial^2}{\partial x^2}\psi_1(x,t)+\left[g_{11}\,|\psi_1(x,t)|^2+g_{12}\,|\psi_2(x,t)|^2\right]\,\psi_1(x,t)&=&0,\label{cnlse1}\\\label{cnlse2}
i\, \frac{\partial}{\partial t}\psi_2(x,t)+g_{20}\,\frac{\partial^2}{\partial x^2}\psi_2(x,t)+\left[g_{21}\,|\psi_1(x,t)|^2+g_{22}\,|\psi_2(x,t)|^2\right]\,\psi_2(x,t)&=&0,
\end{eqnarray}
where, $\psi_1(x,t)$ and $\psi_2(x,t)$ are complex functions, and $g_{10},\,g_{20},\,g_{11},\,g_{12},\,g_{21}$ and $g_{22}$ are real constants.
The two components of the general solution are written in the cartesian complex form
\begin{eqnarray}
\psi_1(x,t)=u_1(x,t)+i\,v_1(x,t),
\end{eqnarray}
\begin{eqnarray}
\psi_2(x,t)=u_2(x,t)+i\,v_2(x,t),
\end{eqnarray}
 where $u_1(x,t)$, $v_1(x,t)$, $u_2(x,t)$, and $v_2(x,t)$ being real functions. Inserting in (\ref{cnlse1}) and (\ref{cnlse2}), generates the following four equations from the real and imaginary parts
\begin{eqnarray}\label{m12}
g_{11}\,u_1^3(x,t)+g_{12}\,u_1(x,t)\,u_2^2(x,t)+g_{11}\,u_1(x,t)\,v_1^2(x,t)+g_{12}\,u_1(x,t)\,v_2^2(x,t)\nonumber&-&\\\frac{\partial}{\partial t}v_1(x,t)+g_{10}\,\frac{\partial^2}{\partial x^2}u_1(x,t)&=&0,\\
g_{11}\,u_1^2(x,t)\,v_1(x,t)+g_{12}\,u_2^2(x,t)\,v_1(x,t)+g_{11}\,v_1^3(x,t)+g_{12}\,v_1(x,t)\,v_2^2(x,t)\nonumber&+&\\\frac{\partial}{\partial t}u_1(x,t)+g_{10}\frac{\partial^2}{\partial x^2}v_1(x,t)&=&0,\\
g_{22}\,u_2^3(x,t)+g_{21}\,u_1^2(x,t)\,u_2(x,t)+g_{21}\,u_2(x,t)\,v_1^2(x,t)+g_{22}\,u_2(x,t)\,v_2^2(x,t)\nonumber&-&\\\frac{\partial}{\partial t}v_2(x,t)+g_{20}\,\frac{\partial^2}{\partial x^2}u_2(x,t)&=&0,\\
g_{21}\,u_1^2(x,t)\,v_2(x,t)+g_{22}\,u_2^2(x,t)\,v_2(x,t)+g_{22}\,v_2^3(x,t)+g_{21}\,v_1^2(x,t)\,v_2(x,t)\nonumber&+&\\\frac{\partial}{\partial t}u_2(x,t)+g_{20}\frac{\partial^2}{\partial x^2}v_2 (x,t)&=&0.\label{m22}
\end{eqnarray}
Then we expand $u_1(x,t)$, $v_1(x,t)$, $u_2(x,t)$, and $v_2(x,t)$ in the following power series
\begin{eqnarray}
u_1(x,t)&=&\sum_{l=0}^{s}a_{l}(x)\,{\Delta t}^{\,l},\label{seru2}\\
v_1(x,t)&=&\sum_{l=0}^{s}b_{l}(x)\,{\Delta t}^{\,l}\label{serv2},\\
u_2(x,t)&=&\sum_{l=0}^{s}c_{l}(x)\,{\Delta t}^{\,l},\label{seru2}\\
v_2(x,t)&=&\sum_{l=0}^{s}d_{l}(x)\,{\Delta t}^{\,l}\label{serv2},
\end{eqnarray}
where the function
coefficients are defined by
\begin{eqnarray}
a_l(x)&=&\frac{1}{l!}\,\frac{\partial^l}{\partial t^l}u_1(x,t)|_{t=0}\label{uvexp11},\\
b_l(x)&=&\frac{1}{l!}\,\frac{\partial^l}{\partial t^l}v_1(x,t)|_{t=0}\label{uvexp22},\\
c_l(x)&=&\frac{1}{l!}\,\frac{\partial^l}{\partial t^l}u_1(x,t)|_{t=0}\label{uvexp11},\\
d_l(x)&=&\frac{1}{l!}\,\frac{\partial^l}{\partial t^l}v_1(x,t)|_{t=0}\label{uvexp22}.
\end{eqnarray}
Substituting these series expansions in Eqs. (\ref{m12})-(\ref{m22}),
recursion relations for the coefficient functions  $a_l(x)$, $b_l(x)$, $c_l(x)$, and $d_l(x)$ are derived and listed in Appendix \ref{appcoupled}  for convenience, as they turn out to be lengthy. The boundary points with  $1\le i\le (p-1)/2$ and $n_x-(p-1)/2< i\le n_x$, are set to the values
\begin{eqnarray}
a_l^i&=&\frac{1}{l!}\frac{\partial^l}{\partial t^l}u_1(x_i,t)|_{t=0}\label{req1b2},\\
b_l^i&=&\frac{1}{l!}\frac{\partial^l}{\partial t^l}v_1(x_i,t)|_{t=0}\label{req2b2},\\
c_l^i&=&\frac{1}{l!}\frac{\partial^l}{\partial t^l}u_2(x_i,t)|_{t=0}\label{req1b22},\\
d_l^i&=&\frac{1}{l!}\frac{\partial^l}{\partial t^l}v_2(x_i,t)|_{t=0}\label{req2b22}.
\end{eqnarray}
The algorithm described in Section \ref{numsec} is then applied to calculate the evolution of the profiles $u_1(x,t)$, $v_1(x,t)$, $u_2(x,t)$, and $v_2(x,t)$.\\\\
\textbf{\textit{Example:  Dark-bright soliton}}\\
The exact dark-bright soliton solution of the two-coupled NLSE, Eqs. (\ref{cnlse1}) and  (\ref{cnlse2}),  is given by
\begin{eqnarray}\label{DB}
\psi_1(x,t)&=&A_0\,\sqrt{\frac{g_{12}\,g_{20}-g_{10}\,g_{22}}{g_{11}\,g_{20}-g_{10}\,g_{21}}}\,\text{tanh}\left\{A_0\sqrt{\frac{g_{12}\,g_{21}-g_{11}\,g_{22}}{2\left(g_{10}\,g_{21}-g_{11}\,g_{20}\right)}}\left[(x-x_0)-k\,(t-t_0)\right]\right\}\\&&\times e^{i\left\{(t-t_0)\left[\frac{g_{11}\,A_0^2\left(g_{12}\,g_{20}-g_{10}\,g_{22}\right)}{g_{11}\,g_{20}-g_{10}\,g_{21}}-\frac{k^2}{4\,g_{10}}\right]+\frac{k\,(x-x_0)}{2\,g_{10}}\right\}},
\nonumber\\
\psi_2(x,t)&=&A_0\,\text{sech}\left\{A_0\sqrt{\frac{g_{12}\,g_{21}-g_{11}\,g_{22}}{2\left(g_{10}\,g_{21}-g_{11}\,g_{20}\right)}}\left[(x-x_0)-k\,(t-t_0)\right]\right\}\\&\times& e^{\frac{i\left\{2g_{12}g_{20}^2g_{21}A_0^2(t-t_0)+2g_{11}g_{20}g_{20}g_{22}A_0^2(t-t_0)+kg_{11}g_{20}\left[2(x-x_0)-k(t-t_0)\right]+4g_{10}g_{21}g_{20}g_{22}A_0^2(t-t_0)+kg_{10}g_{21}\left[2(x-x_0)+k(t-t_0)\right]\right\}}{4g_{20}\left(g_{11}g_{20}-g_{10}g_{21}\right)}},\nonumber\\
\end{eqnarray}
where $A_0$, $k$, $x_0$, and $t_0$ are arbitrary real constants and all quantities under the square root must be positive.

\begin{figure}
    \includegraphics[scale=0.60]{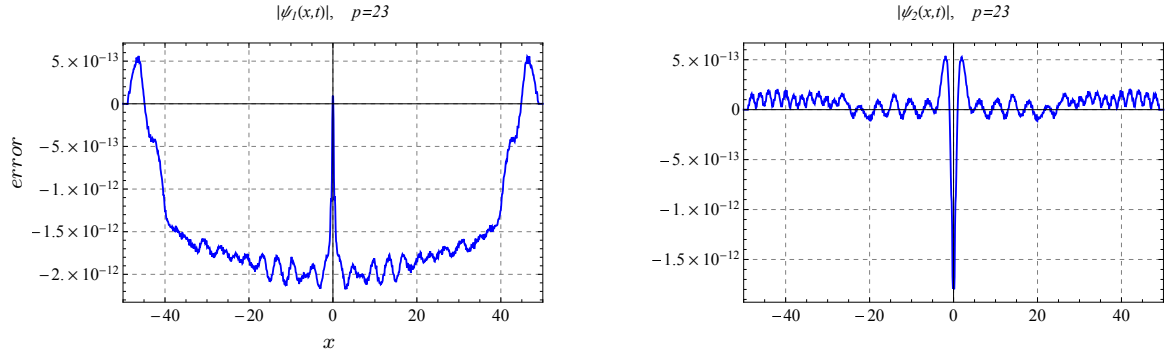}
    \begin{picture}(5,5)(5,5)
    \put(-434,55) {\rotatebox{90}{\scriptsize$error$}}
    \end{picture}
    \begin{picture}(5,5)(5,5)
    \put(-325,-5) {\scriptsize$x$}
    \end{picture}\vspace{0.2cm}
    \caption{Error in the numerical solution corresponding to the dark-bright soliton of the two-coupled NLSE (\ref{cnlse1}) and (\ref{cnlse2}) at $t=40$ using $p=23$ and $s=4$. Left:
    Error associated with $\psi_1(x,t)$. Right:
    Error associated with $\psi_2(x,t)$. Parameters:
    $n_x=1000$, $n_t=80000$, $L=100$, $\Delta t=0.0005$,
    $\Delta x=L/(n_x-1)$, $k=0.0$, $x_0=t_0=0.0$, $A_0=1.0$, $g_{10}=1/2$, $g_{11}=-1$, $g_{12}=1/2$, $g_{20}=1/2$, $g_{21}=-1/2$, $g_{22}=1$. }
    \label{fig20}
\end{figure}

In Fig. \ref{fig20}, we plot the error versus $x$ in the numerical solution of both components, $\psi_1(x,t)$ and $\psi_2(x,t)$, at $t=40$ and using $p=23$ and $s=4$. Clearly, the error in both components is extremely small even for such a long evolution time.

\section{Conclusions and Outlook}
\label{concsec}
We have presented a high accuracy numerical method that solves the initial value problem of the fundamental, inhomogeneous, and coupled NLSE. The method employed an iterative power series for time stepping and a multi-point formula for the spacial discretisation of the second derivative. The method is characterised by a systematic increase in accuracy in terms of the two parameters $s$ and $p$, representing the maximum power of the time power series and number of points in the multi-point formula, respectively. As a result, the accuracy was shown via some examples to reach the machine precision in a rather short computing time. 

Detailed analysis of the different sources of error was performed. Errors arising from the time power series, $error_s$, and the multi-point formula, $error_p$, were characterised and accounted theoretically. Error from the boundaries, which is more significant for solutions with uniform background, was shown to reduce with higher order of the time power series and larger system size. 

We compared the present method with two other methods representing the finite difference and spectral methods. For finite difference methods, the G-FDTD method was selected for comparison due to its high accuracy and similarity to our method. For the spectral methods, we compared our results with the Fourier split-step method.  We have pointed out the similarities and differences between our method and the G-FDTD method and shown that the present method is characterised by a faster computing time and higher convergence rate. The accuracy and computing time exceed by far those of the split-step method. The method is extended to the inhomogeneous NLSE and applied to the scattering of a bright soliton by a reflectionless potential. This example showed the importance of high accuracy to capture the correct scattering outcome near the critical value of soliton speed separating quantum reflection from transmission. We have also generalised the method to the two-coupled NLSE and considered the example of  dark-bright soliton. 

It is straightforward to extend the method to higher dimensions. The method can also be extended to NLSE with higher order terms such as third dispersion and  Raman scattering, etc. It should be noted however, that the method does not apply to time-dependent potentials and time-dependent coefficients. In obtaining the recursion relations of the time power series, this was implicitly assumed. Extending the method to evolution equations with higher time derivatives is also possible, but with a different recursive structure.

In conclusion, we believe the method presented here will be very useful for realistic efficient numerical solutions of nonlinear evolution equations.

\section*{Acknowledgment} The authors acknowledge the support of UAE University through grants UAEU-UPAR (1) 2019 and UAEU-UPAR (11) 2019.

\clearpage
\section*{Tables}
\begin{table}[H]
	\centering
	\begin{math}
	\begin{array}{|cccc|}
	\hline
	&&\Delta t=1\times10^{-4}&\\\hline
	&&{\rm Present \,\,method}&\\
	\hline
	n_x & {\rm error} &{\rm CPU\, time (s)} & {\rm convergence\, rate}  \\
	\hline
	100 & 2.04764\times 10^{-2} & 0.193&  \\
	200 & 1.24616\times 10^{-3} & 0.332& 4.06753 \\
	300 & 2.49757\times 10^{-4} & 0.465& 3.98046 \\
	400 & 7.95110\times 10^{-5} & 0.595& 3.9902 \\
	500 & 3.26624\times 10^{-5} & 1.030 &
	3.99593 \\
	600 & 1.57757\times 10^{-5} & 0.990& 3.99889 \\
	700 & 8.5228\times 10^{-6} & 1.072 &
	4.00048 \\
	800 & 4.9986\times 10^{-6} & 1.186 &
	4.00136 \\
	900 & 3.1216\times 10^{-6} & 1.463 &
	4.00185 \\
	1000 & 2.0486\times 10^{-6} & 1.707 &
	4.00212 \\
	\hline
	&&{\rm G-FDTD}&\\
	\hline
	n_x & {\rm error} & {\rm CPU\, time (s)} &{\rm  convergence\,\, rate}  \\
	100 & 2.11658\times 10^{-2} & 0.415&  \\
	200 & 1.25608\times 10^{-3} & 0.796& 4.07473 \\
	300 & 2.68115\times 10^{-4} & 1.195 &
	3.80880 \\
	400 & 1.43054\times 10^{-4} & 1.572 &
	2.18364 \\
	\hline
	\end{array}
	\end{math}
	\caption{Error data for bright soliton given by (\ref{bs}) with $g_1=-1,\,g_2=-2, A_0=1,\, k=4,\,x_0=-10$, $p=5,\,s=3$. The G-FDTD code diverges for $n_x>400$.}
	\label{table1}
\end{table}

\begin{table}[H]
		\centering
	\begin{math}
	\begin{array}{|cccc|}
	\hline
	&&\Delta t=1\times10^{-6}&\\\hline
	&&{\rm Present\,\, method}&\\
	\hline
	n_x & {\rm error} &{\rm CPU\, time (s)} & {\rm convergence\, rate}  \\
	\hline
	100 & 2.04764\times 10^{-2} & 1.600\times
	10^1 &  \\
	200 & 1.24616\times 10^{-3} & 3.387\times
	10^1 & 4.06753 \\
	300 & 2.49757\times 10^{-4} & 4.998\times
	10^1 & 3.98046 \\
	400 & 7.95110\times 10^{-5} & 6.462\times
	10^1 & 3.9902 \\
	500 & 3.26623\times 10^{-5} & 9.762\times
	10^1 & 3.99593 \\
	600 & 1.57757\times 10^{-5} & 1.015\times
	10^2 & 3.99889 \\
	700 & 8.5228\times 10^{-6} & 1.167\times
	10^2 & 4.00048 \\
	800 & 4.9986\times 10^{-6} & 1.321\times
	10^2 & 4.00136 \\
	900 & 3.1216\times 10^{-6} & 1.491\times
	10^2 & 4.00185 \\
	1000 & 2.0485\times 10^{-6} & 1.833\times
	10^2 & 4.00213 \\
	\hline
	&&{\rm G-FDTD}&\\
	\hline
	n_x & {\rm error} & {\rm CPU\, time (s)} &{\rm  convergence\,\, rate}  \\
	100 & 2.11482\times 10^{-2} & 4.066\times
	10^1 &  \\
	200 & 1.26552\times 10^{-3} & 8.348\times
	10^1 & 4.06273 \\
	300 & 2.52158\times 10^{-4} & 1.274\times
	10^2 & 3.9786 \\
	400 & 7.99569\times 10^{-5} & 1.707\times
	10^2 & 3.99249 \\
	500 & 3.27135\times 10^{-5} & 2.095\times
	10^2 & 4.00504 \\
	600 & 1.57349\times 10^{-5} & 2.518\times
	10^2 & 4.01439 \\
	700 & 8.4839\times 10^{-6} & 2.901\times
	10^2 & 4.00714 \\
	800 & 5.0105\times 10^{-6} & 3.305\times
	10^2 & 3.94393 \\
	900 & 3.2231\times 10^{-6} & 3.769\times
	10^2 & 3.7457 \\
	1000 & 2.2738\times 10^{-6} & 4.246\times
	10^2 & 3.31165 \\
	\hline
	\end{array}
	\end{math}
	\caption{Error data for bright soliton given by (\ref{bs}) with $g_1=-1,\,g_2=-2, A_0=1,\, k=4,\,x_0=-10$, $p=5,\,s=3$. }
	\label{table2}
\end{table}

\begin{table}[H]
		\centering
	\begin{math}
	\begin{array}{|cccc|}
	\hline
	&&p=11&\\
	\hline
	n_x & {\rm error} &{\rm CPU\, time\, (s)} & {\rm convergence\, rate}  \\
	\hline
	50 & 7.41644\times 10^{-2} & 0.246&  \\
	100 & 8.4102\times 10^{-4} & 0.384& 6.55562 \\
	150 & 2.38048\times 10^{-5} & 0.544& 8.86396 \\
	200 & 1.67858\times 10^{-6} & 0.679& 9.27171 \\
	250 & 2.0163\times 10^{-7} & 0.777& 9.53988 \\
	300 & 3.46857\times 10^{-8} & 0.888& 9.68916 \\
	350 & 7.7177\times 10^{-9} & 0.925& 9.77911 \\
	400 & 2.0819\times 10^{-9} & 1.009 &
	9.83856 \\
	450 & 6.517\times 10^{-10} & 1.136 &
	9.88372 \\
	500 & 2.296\times 10^{-10} & 1.319 &
	9.92345 \\
	\hline
	&&p=15&\\\hline
	50 & 6.14867\times 10^{-2} & 0.381&  \\
	100 & 2.5016\times 10^{-4} & 0.428& 8.05579 \\
	150 & 2.63863\times 10^{-6} & 0.711& 11.3185 \\
	200 & 8.60224\times 10^{-8} & 0.754& 11.9689 \\
	250 & 5.29238\times 10^{-9} & 0.946 & 12.5517 \\
	300 & 5.04834\times 10^{-10} & 1.069 &
	12.9353 \\
	350 & 6.64545\times 10^{-11} & 1.138 &
	13.1947 \\
	400 & 1.11947\times 10^{-11} & 1.321 &
	13.3739 \\
	450 & 2.29149\times 10^{-12} & 1.474 &
	13.4992 \\
	500 & 5.48892\times 10^{-13} & 1.747 &
	13.5921 \\
	\hline
	&&p=23&\\\hline
	50 & 4.89514\times 10^{-2} & 0.517&  \\
	100 & 8.71723\times 10^{-5} & 0.695 & 9.26496 \\
	150 & 3.00017\times 10^{-7} & 0.894 & 14.1034 \\
	200 & 4.07153\times 10^{-9} & 0.950& 15.0330 \\
	250 & 1.02235\times 10^{-10} & 1.160 &
	16.5858 \\
	300 & 4.04200\times 10^{-12} & 1.204 &
	17.7837 \\
	350 & 2.30026\times 10^{-13} & 1.459 &
	18.6516 \\
	400 & 1.77624\times 10^{-14} & 1.737 &
	19.2312 \\
	450 & 3.50575\times 10^{-15} & 2.046 &
	13.8093 \\
	500 & 2.67075\times 10^{-15} & 2.280 &
	2.58748 \\
	\hline
	\end{array}
	\end{math}
	\caption{Error data of the present method for bright soliton given by (\ref{bs}) with $g_1=-1,\,g_2=-2, A_0=1,\, k=4,\,x_0=-10$, $\Delta t=10^{-4}$, $s=4$. }
	\label{table3}
\end{table}

\begin{table}[H]
		\centering
	\begin{math}
	\begin{array}{|cccc|}
	\hline
	&&{\rm Present \,\,method}&\\\hline
	n_x & {\rm error} &{\rm CPU\, time\, (s)} & {\rm convergence\, rate}  \\
	\hline
	100 & 7.00299\times 10^{-3} & 3.\times 10^{-2} &  \\
	200 & 3.01673\times 10^{-4} & 3.6\times 10^{-2} & 4.56963 \\
	300 & 6.11528\times 10^{-5} & 5.7\times 10^{-2} & 3.95232 \\
	400 & 1.96238\times 10^{-5} & 7.\times 10^{-2} & 3.96245 \\
	500 & 8.0955\times 10^{-6} & 1.07\times 10^{-1} & 3.97687 \\
	600 & 3.9212\times 10^{-6} & 1.07\times 10^{-1} & 3.98333 \\
	700 & 2.1266\times 10^{-6} & 1.4\times 10^{-1} & 3.97537 \\
	800 & 1.2586\times 10^{-6} & 1.25\times 10^{-1} & 3.93352 \\
	900 & 8.032\times 10^{-7} & 1.4\times 10^{-1} & 3.81766 \\
	1000 & 5.518\times 10^{-7} & 1.77\times 10^{-1} & 3.56716 \\
	\hline
	&&{\rm G-FDTD}&\\
	\hline 100 & 7.37839\times 10^{-3} & 2.3\times 10^{-2} &  \\
	200 & 4.71382\times 10^{-4} & 3.8\times 10^{-2} & 3.96834 \\
	300 & 2.4642\times 10^{-4} & 5.6\times 10^{-2} & 1.59972 \\
	400 & 2.11859\times 10^{-4} & 6.6\times 10^{-2} & 5.2529\times
	10^{-1} \\
	500 & 2.02573\times 10^{-4} & 8.0\times 10^{-2} & 2.0086\times
	10^{-1} \\
	600 & 1.99137\times 10^{-4} & 1.13\times 10^{-1} & 9.385\times
	10^{-2} \\
	700 & 1.97577\times 10^{-4} & 1.11\times 10^{-1} & 5.100\times
	10^{-2} \\
	800 & 1.9676\times 10^{-4} & 1.25\times 10^{-1} & 3.105\times
	10^{-2} \\
	900 & 1.96282\times 10^{-4} & 1.5\times 10^{-1} & 2.063\times
	10^{-2} \\
	1000 & 1.95979\times 10^{-4} & 1.6\times 10^{-1} & 1.469\times
	10^{-2} \\
	\hline
	\end{array}
	\end{math}
	\caption{Error data for dark soliton given by (\ref{dark}) with $g_1=1/2,\,g_2=-4, A_0=1,\, k=4,\,x_0=-10$, $p=5,\,s=3$. Note that in Ref. \cite{meth3}, the function given by (24) is not a dark soliton solution to Eq. (\ref{nlse}) (also Eq. (1) in Ref. \cite{meth3}) with the given values of $g_1=1$ and $g_2=-4$. Instead, $g_1=1/2,\,g_2=-4$ are the correct ones.  }
	\label{table4}
\end{table}

\clearpage
\section*{Figures}



\clearpage

\appendix
\section{Deriving the $p$-point central difference formula of the second derivative}
\label{appa}
We first derive the $p$-point central difference formula
that approximates a second derivative for any odd integer $p$. Then we present specific
examples of $p=5,\dots,23$.

Consider the Taylor expansion
\begin{equation}
f(x+j\Delta
x)=\sum_{i=0}^{p}\frac{(j\Delta x)^i}{i!}\frac{d^if(x)}{dx^i}+\mathcal{O}\left(\frac{(j\Delta x)^{p+1}}{(p+1)!}\frac{d^{p+1}f(x)}{dx^{p+1}}\right)\label{eqapp1},
\end{equation}
and the linear combination
\begin{equation}
\sum_{j=1}^{(p-1)/2}C_j\left[f(x+j\Delta x)+f(x-j\Delta x)\right],
\end{equation}
where $C_j$ are arbitrary real coefficients. Substituting the Taylor
expansion (\ref{eqapp1}) in the last expression and then exchanging the summations, it takes the form
\begin{eqnarray}
&&\sum_{j=1}^{(p-1)/2}C_j\left[f(x+j\Delta x)+f(x-j\Delta
x)\right]\nonumber\\&=&2\sum_{i=0}^{(p-1)/2}\left[\sum_{j=1}^{(p-1)/2}C_j\frac{(j\Delta
x)^{2i}}{(2i)!}\right]\frac{df^{2i}(x)}{dx^{2i}}+\mathcal{O}\left(\frac{(j\Delta x)^{p+1}}{(p+1)!}\frac{d^{p+1}f(x)}{dx^{p+1}}\right)\label{eqapp4}.
\end{eqnarray}
Requesting all coefficients of $d^{2i}f(x)/dx^{2i}$ for $i>1$ to
vanish, gives the following system of equations
\begin{equation}
\sum_{j=1}^{(p-1)/2}C_j\frac{(j\Delta x)^{2i}}{(2i)!}=0\label{eqapp2}
\end{equation}
for $i=2,3,\dots,(p-1)/2$. Solving (\ref{eqapp4}) for the remaining second
derivative, we obtain
\begin{eqnarray}
f_p^{\prime\prime}&\equiv&\frac{d^2f(x)}{dx^2}=\frac{1}{\Delta
x^2\sum_{j=1}^{(p-1)/2}C_j{j^2}}\nonumber\\&\times&\sum_{j=1}^{(p-1)/2}C_j\left[f(x+j\Delta
x)+f(x-j\Delta
x)-2f(x)\right]\nonumber\\&+&\mathcal{O}\left(\frac{(j\Delta x)^{p+1}}{(p+1)!}\frac{d^{p+1}f(x)}{dx^{p+1}}\right)\label{eqapp3},
\end{eqnarray}
where we denote the $p$-point second derivative as $f^{\prime\prime}_p$. The last term gives the order of  error in the approximation.
As an example, we take the $p=5$ case. The system (\ref{eqapp2})
reduces to a single equation for $i=2$, namely $c_1+16c_2=0$. Solving
for $c_2$ and substituting back in (\ref{eqapp3}),
gives the 5-point central formula as listed below.

For convenience, we calculate the $p$-point formulae up to $p=23$:
\begin{eqnarray}
f^{\prime\prime}_5(x)&=&\frac{-f(x+2\,\Delta x)+16\,f(x+\Delta
x)-30\,f(x)+16\,f(x-\Delta x)-f(x-2\,\Delta x)}{{12\,\Delta x}^2},
\end{eqnarray}
\begin{eqnarray}
f^{\prime\prime}_7(x)&=&\frac{1}{{180\,\Delta
x}^2}\Big[2\,f(x+3\,\Delta x)-27\,f(x+2\,\Delta x)+270\,f(x+\Delta
x)-490\,f(x)\\\nonumber&&+270\,f(x-\Delta x)-27\,f(x-2\,\Delta
x)+2\,f(x-3\,\Delta x)\Big],
\end{eqnarray}
\begin{eqnarray}
f^{\prime\prime}_9(x)&=&\frac{1}{{5040\,\Delta
x}^2}\Big[-9\,f(x+4\,\Delta x)+128\,f(x+3\,\Delta
x)-1008\,f(x+2\,\Delta x)\\\nonumber&&+8064\,f(x+\Delta
x)-14350\,f(x)+8064\,f(x-\Delta x)-1008\,f(x-2\,\Delta
x)\\\nonumber&&+128\,f(x-3\,\Delta x)-9\,f(x-4\,\Delta x)\Big],
\end{eqnarray}
\begin{eqnarray}
f^{\prime\prime}_{11}(x)&=&\frac{1}{{25200\,\Delta
x}^2}\Big[8\,f(x+5\,\Delta x)-125\,f(x+4\,\Delta
x)+1000\,f(x+3\,\Delta x)\\\nonumber&&-6000\,f(x+2\,\Delta
x)+42000\,f(x+\Delta x)-73766\,f(x)+42000\,f(x-\Delta
x)\\\nonumber&&-6000\,f(x-2\,\Delta x)+1000\,f(x-3\,\Delta
x)-125\,f(x-4\,\Delta x)+8\,f(x-5\,\Delta x)\Big],
\end{eqnarray}
\begin{eqnarray}
f^{\prime\prime}_{13}(x)&=&\frac{1}{{831600\,\Delta
x}^2}\Big[-50\,f(x+6\,\Delta x)+864\,f(x+5\,\Delta
x)-7425\,f(x+4\,\Delta x)\\\nonumber&&+44000\,f(x+3\,\Delta
x)-222750\,f(x+2\,\Delta x)+1425600\,f(x+\Delta
x)-2480478\,f(x)\\\nonumber&&+1425600\,f(x-\Delta
x)-222750\,f(x-2\,\Delta x)+44000\,f(x-3\,\Delta
x)\\\nonumber&&-7425\,f(x-4\,\Delta x)+864\,f(x-5\,\Delta
x)-50\,f(x-6\,\Delta x)\Big],
\end{eqnarray}
\begin{eqnarray}
f^{\prime\prime}_{15}(x)&=&\frac{1}{{75675600\,\Delta
x}^2}\Big[900\,f(x+7\,\Delta x)-17150\,f(x+6\,\Delta
x)+160524\,f(x+5\,\Delta x)\\\nonumber&&-1003275\,f(x+4\,\Delta
x)+4904900\,f(x+3\,\Delta x)-22072050\,f(x+2\,\Delta
x)\\\nonumber&&+132432300\,f(x+\Delta
x)-228812298\,f(x)+132432300\,f(x-\Delta
x)\\\nonumber&&-22072050\,f(x-2\,\Delta x)+4904900\,f(x-3\,\Delta
x)-1003275\,f(x-4\,\Delta x)\\\nonumber&&+160524\,f(x-5\,\Delta
x)-17150\,f(x-6\,\Delta x)+900\,f(x-7\,\Delta x)\Big],
\end{eqnarray}
\begin{eqnarray}
f^{\prime\prime}_{17}(x)&=&\frac{1}{{302702400\,\Delta
x}^2}\Big[-735\,f(x+8\,\Delta x)+15360\,f(x+7\,\Delta
x)\\\nonumber&&-156800\,f(x+6\,\Delta x)+1053696\,f(x+5\,\Delta
x)-5350800\,f(x+4\,\Delta x)\\\nonumber&&+22830080\,f(x+3\,\Delta
x)-94174080\,f(x+2\,\Delta x)+538137600\,f(x+\Delta
x)\\\nonumber&&-924708642\,f(x)+538137600\,f(x-\Delta
x)-94174080\,f(x-2\,\Delta x)\\\nonumber&&+22830080\,f(x-3\,\Delta
x)-5350800\,f(x-4\,\Delta x)+1053696\,f(x-5\,\Delta
x)\\\nonumber&&-156800\,f(x-6\,\Delta x)+15360\,f(x-7\,\Delta
x)-735\,f(x-8\,\Delta x)\Big],
\end{eqnarray}
\begin{eqnarray}
f^{\prime\prime}_{19}(x)&=&\frac{1}{{15437822400\,\Delta
x}^2}\Big[7840\,f(x+9\,\Delta x)-178605\,f(x+8\,\Delta
x)\\\nonumber&&+1982880\,f(x+7\,\Delta x)-14394240\,f(x+6\,\Delta
x)+77728896\,f(x+5\,\Delta x)\\\nonumber&&-340063920\,f(x+4\,\Delta
x)+1309875840\,f(x+3\,\Delta x)-5052378240\,f(x+2\,\Delta
x)\\\nonumber&&+27788080320\,f(x+\Delta
x)-47541321542\,f(x)+27788080320\,f(x-\Delta
x)\\\nonumber&&-5052378240\,f(x-2\,\Delta
x)+1309875840\,f(x-3\,\Delta x)-340063920\,f(x-4\,\Delta
x)\\\nonumber&&+77728896\,f(x-5\,\Delta x)-14394240\,f(x-6\,\Delta
x)+1982880\,f(x-7\,\Delta x)\\\nonumber&&-178605\,f(x-8\,\Delta
x)+7840\,f(x+9\,\Delta x)\Big],
\end{eqnarray}
\begin{eqnarray}
f^{\prime\prime}_{21}(x)&=&\frac{1}{{293318625600\,\Delta
x}^2}\Big[-31752\,f(x+10\,\Delta x)+784000\,f(x+9\,\Delta
x)\\\nonumber&&-9426375\,f(x+8\,\Delta x)+73872000\,f(x+7\,\Delta
x)-427329000\,f(x+6\,\Delta
x)\\\nonumber&&+1969132032\,f(x+5\,\Delta
x)-7691922000\,f(x+4\,\Delta x)+27349056000\,f(x+3\,\Delta
x)\\\nonumber&&-99994986000\,f(x+2\,\Delta
x)+533306592000\,f(x+\Delta
x)-909151481810\,f(x)\\\nonumber&&+533306592000\,f(x-\Delta
x)-99994986000\,f(x-2\,\Delta x)+27349056000\,f(x-3\,\Delta
x)\\\nonumber&&-7691922000\,f(x-4\,\Delta
x)+1969132032\,f(x-5\,\Delta x)-427329000\,f(x-6\,\Delta
x)\\\nonumber&&+73872000\,f(x-7\,\Delta x)-9426375\,f(x-8\,\Delta
x)+784000\,f(x+9\,\Delta x)\\\nonumber&&-31752\,f(x-10\,\Delta
x)\Big],
\end{eqnarray}
\begin{eqnarray}
f^{\prime\prime}_{23}(x)&=&\frac{1}{{3226504881600\,\Delta
x}^2}\Big[75600\,f(x+11\,\Delta x)-2012472\,f(x+10\,\Delta
x)\\\nonumber&&+26087600\,f(x+9\,\Delta x)-220114125\,f(x+8\,\Delta
x)+1365606000\,f(x+7\,\Delta
x)\\\nonumber&&-6691469400\,f(x+6\,\Delta
x)+27301195152\,f(x+5\,\Delta x)-97504268400\,f(x+4\,\Delta
x)\\\nonumber&&+325014228000\,f(x+3\,\Delta
x)-1137549798000\,f(x+2\,\Delta x)+5915258949600\,f(x+\Delta
x)\\\nonumber&&-10053996959110\,f(x)+5915258949600\,f(x-\Delta
x)-1137549798000\,f(x-2\,\Delta
x)\\\nonumber&&+325014228000\,f(x-3\,\Delta
x)-97504268400\,f(x-4\,\Delta x)+27301195152\,f(x-5\,\Delta
x)\\\nonumber&&-6691469400\,f(x-6\,\Delta
x)+1365606000\,f(x-7\,\Delta x)-220114125\,f(x-8\,\Delta
x)\\\nonumber&&+26087600\,f(x+9\,\Delta x)-2012472\,f(x-10\,\Delta
x)+75600\,f(x-11\,\Delta x)\Big].
\end{eqnarray}

\section{Recursion relations of the two-coupled NLSE}
\label{appcoupled}
Recursion relations for the  two-coupled NLSE (\ref{m12}-\ref{m22}) up to the fourth order, $s=4$:
\begin{eqnarray}
a_1&=&-g_{11}\,b_0^3-b_0\left[g_{11}\,a_0^2+g_{12}\left(c_0^2+\,d_0^2\right)\right]-g_{10}\,b_0^{\prime\prime},\\
a_2&=&\frac{1}{2}\Big[-2g_{11}\,a_0\,a_1b_0-g_{11}\,a_0^2\,b_1-3g_{11}\,b_0^2\,b_1-g_{12}\,b_1\,c_0^2-g_{12}\,b_1d_0^2\nonumber\\&&-2g_{12}\,b_0\left(c_0\,c_1+d_0\,d_1\right)-g_{10}b_1^{\prime\prime}\Big],\\
a_3&=&\frac{1}{3}\Big(-g_{11}\,a_1^2\,b_0-2\,g_{11}\,a_0\,a_2\,b_0-2\,g_{11}\,a_0\,a_1\,b_1-3\,g_{11}\,b_0\,b_1^2-g_{11}\,a_0^2\,b_2-3\,g_{11}\,b_0^2\,b_2\nonumber\\&&-g_{12}\,b_2\,c_0^2-2\,g_{12}\,b_1\,c_0\,c_1-g_{12}\,b_0\,c_1^2-2\,g_{12}\,b_0\,c_0\,c_2-g_{12}\,b_2\,d_0^2-2\,g_12\,b1\,d_0\,d_1\nonumber\\&&-g_{12}\,b_0\,d_1^2-2\,g_{12}\,b_0\,d_0\,d_2-g_{10}\,b_2^{\prime\prime}\Big),\\
a_4&=&\frac{1}{4}\Big[-g_{11}\,b_1^3-2\,g_{11}\,a_0\left(a_3\,b_0+a_2\,b_1\right)-2\,g_{11}\,a_1\left(a_2\,b_0+a_0\,b_2\right)-g_{11}\,a_0^2\,b_3-3\,g_{11}\,b_0^2\,b_3\nonumber\\&&-g_{12}\,b_3\,c_0^2-2\,g_{12}\,b_2\,c_0\,c_1-g_{12}\,b_1\,c_1^2-2\,g_{12}\,b_0\,c_1\,c_2-2\,g_{12}\,b_0\,c_0\,c_3-g_{12}\,b_3\,d_0^2\nonumber\\\nonumber&&-2\,g_{12}\,b_2\,d_0\,d_1-2\,g_{12}\,b_0\,d_1\,d_2-2\,g_{12}\,b_0\,d_0\,d_3-g_{10}\,b_3^{\prime\prime}\nonumber\\&&-b_1\left(g_{11}\,a_1^2+6\,g_{11}\,b_0\,b_2+2\,g_{12}\,c_0\,c_2+g_{12}\,d_1^2+2\,g_{12}\,d_0\,d_2\right)\Big]\\
b_1&=&g_{11}\,a_0^3+a_0\left[g_{11}\,b_0^2+g_{12}\left(c_0^2+\,d_0^2\right)\right]+g_{10}\,a_0^{\prime\prime},\\
b_2&=&\frac{1}{2}\Big\{3\,g_{11}\,a_0^2\,a_1+a_1\left[g_{11}\,b_0^2+g_{12}\left(c_0^2+\,d_0^2\right)\right]+2\,a_0\left(g_{11}\,b_0+g_{12}\,c_0\,c_1+g_{12}\,d_0\,d_1\right)\nonumber\\&&+g_{10}\,a_1^{\prime\prime}\Big\},\\b_3&=&\frac{1}{3}\Big\{3g_{11}\,a_0^2+2g_{11}\,a_1\,b_0\,b_1+2g_{12}\,a_1\,c_0\,c_1+a_2\left[g_{11}\,b_0^2+g_{12}\left(c_0^2+\,d_0^2\right)\right]+2g_{12}\,a_1\,d_0\,d_1\nonumber\\&&+a_0\left(3\,g_{11}\,a_1^2+g_{11}\,b_1^2+2\,g_{11}\,b_0\,b_2+g_{12}\,c_1^2+2\,g_{12}\,c_0\,c_2+g_{12}\,d_0\,d_2\right)+g_{10}\,a_2^{\prime\prime}\Big\},\\b_4&=&\frac{1}{4}\Big\{g_{11}\,a_1^3+3\,g_{11}\,a_0^2\,a_3+g_{11}\,a_3\,b_0^2+2\,g_{11}\,a_2\,b_0\,b_1+g_{12}\,a_3\,c_0^2\nonumber\\&&+2\,g_{12}\,a_2\,c_0\,c_1+g_{12}\,a_3\,d_0^2+2\,g_{12}\,a_2\,d_0\,d_1+g_{10}\,a_3^{\prime\prime}\nonumber\\&&+a_1\left(6\,g_{11}\,a_0\,a_2+g_{11}\,b_1^2+2\,g_{11}\,b_0\,b_2+g_{12}\,c_1^2+2\,g_{12}\,c_0\,c_2+g_{12}\,d_1^2+2\,g_{12}\,d_0\,d_2\right)\nonumber\\&&+2\,a_0\left[g_{11}\,b_1\,b_2+g_{11}\,b_0\,b_3+g_{12}\left(c_1\,c_2+c_0\,c_3+d_1\,d_2+d_0\,d_3\right)\right]\Big\},\\
c_1&=&-g_{21}\,a_0^2\,d_0-g_{21}\,b_0^2\,d_0-g_{22}\,c_0^2\,d_0-g_{22}\,d_0^3-g_{20}\,d_0^{\prime\prime},\\c_2&=&\frac{1}{2}\Big(-2\,g_{21}\,a_0\,a_1\,d_0-2\,g_{21}\,b_0\,b_{1}\,d_0-2\,g_{22}\,c_0\,c_1\,d_0-g_{21}\,a_0^2\,d_1-g_{21}\,b_0^2\,d_1-g_{22}\,c_0^2\,d_1\nonumber\\&&-3\,g_{22}\,d_0^2\,d_1-g_{20}\,d_1^{\prime\prime}\Big),\\
c_3&=&\frac{1}{3}\Big(-g_{21}\,a_1^2\,d_0-2\,g_{21}\,a_0\,a_2\,d_0-g_{21}\,b_1^2\,d_0-2\,g_{21}\,b_0\,b_2\,d_0-g_{22}\,c_1^2\,d_0-2\,g_{22}\,c_0\,c_2\,d_0\nonumber\\&&-2\,g_{21}\,a_0\,a_1\,d_1-2\,g_21\,b_0\,b_1\,d_1-2\,g_{22}\,c_0\,c_1\,d_1-3\,g_{22}\,d_0\,d_1^2-g_{21}\,a_0^2\,d_2-g_{21}\,b_0^2\,d_2\nonumber\\&&-g_{22}\,c_0^2\,d_2-3\,g_{22}\,d_0^2\,d_2-g_{20}\,d_2^{\prime\prime}\Big),\\
c_4&=&\frac{1}{4}\Big[-2g_{21}\,b_1\,b_2\,d_0-2g_{21}\,b_0\,b_3\,d_0-2g_{22}\,c_1\,c_2\,d_0-2g_{22}\,c_0\,c_3\,d_0-g_{21}\,a_1^2\,d_1-g_{21}\,b_1^2\,d_1\nonumber\\&&-2\,g_{21}\,b_0\,b_2\,d_1-g_{22}\,c_1^2\,d_1-2\,g_{22}\,c_0\,c_2\,d_1-g_{22}\,d_1^3-2\,g_{21}\,a_0\left(a_3\,d_0+a_2\,d_1\right)\nonumber\\&&-2\,g_{21}\,b_0\,b_1\,d_2-2\,g_{22}\,c_0\,c_1\,d_2-6\,g_{22}\,d_0\,d_1\,d_2-2\,g_{21}\,a_1\left(a_2\,d_0+a_0\,d_2\right)-g_{21}\,a_0^2\,d_3\nonumber\\&&-g_{21}\,b_0^2\,d_3-g_{22}\,c_0^2\,d_3-3\,d_0^2\,d_3-g_{20}\,d_3^{\prime\prime}\Big],\\
d_1&=&g_{21}\,a_0^2\,c_0+g_{21}\,b_0^2\,c_0+g_{22}\,c_0\,d_0^2+g_{22}\,c_0^3+g_{20}\,c_0^{\prime\prime},\\d_2&=&\frac{1}{2}\Big(2\,g_{21}\,a_0\,a_1\,c_0+2\,g_{21}\,b_0\,b_{1}\,c_0+2\,g_{22}\,c_0\,d_0\,d_1+g_{21}\,a_0^2\,c_1+g_{21}\,b_0^2\,c_1+g_{22}\,d_0^2\,c_1\nonumber\\&&+3\,g_{22}\,c_0^2\,c_1+g_{20}\,c_1^{\prime\prime}\Big),\\
d_3&=&\frac{1}{3}\Big(g_{21}\,a_1^2\,c_0+2\,g_{21}\,a_0\,a_2\,c_0+g_{21}\,b_1^2\,c_0+2\,g_{21}\,b_0\,b_2\,c_0+2\,g_{21}\,a_0\,a_1\,c_1+2\,g_{21}\,b_0\,b_1\,c_1\nonumber\\&&+3\,g_{22}\,c_0\,c_1^2+g_{21}\,a_0^2\,c_2+g_{21}\,b_0^2\,c_2+3\,g_{22}\,c_0^2\,c_2+g_{22}\,c_2\,d_0^2+2\,g_{22}\,c_1\,d_0\,d_1\nonumber\\&&+g_{22}\,c_0\,d_1^2+2\,g_{22}\,c_0\,d_0\,d_2+g_{20}\,c_2^{\prime\prime}\Big),\\
d_4&=&\frac{1}{4}\Big[2g_{21}\,b_1\,b_2\,c_0+2g_{21}\,b_0\,b_3\,c_0+g_{21}\,a_1^2\,c_1+g_{21}\,b_1^2\,c_1+2\,g_{21}\,b_0\,b_2\,c_1+g_{22}\,c_1^3\nonumber\\&&+2\,g_{21}\,a_0\left(a_3\,c_0+a_2\,c_1\right)+2\,g_{21}\,b_0\,b_1\,c_2+6\,g_{22}\,c_0\,c_1\,c_2+2\,g_{21}\,a_1\left(a_2\,c_0+a_0\,c_2\right)\nonumber\\&&+g_{21}\,a_0^2\,c_3+g_{21}\,b_0^2\,c_3+3\,g_{22}\,c_0^2\,c_3+g_{22}\,c_3\,d_0^2+2\,g_{22}\,c_2\,d_0\,d_1+g_{22}\,c_1\,d_1^2\nonumber\\&&+2\,g_{22}\,c_1\,d_0\,d_2+2\,g_{22}\,c_0\,d_1\,d_2+2\,g_{22}\,c_0\,d_0\,d_3+g_{20}\,c_3^{\prime\prime}\Big].
\end{eqnarray}

\end{document}